\documentclass[
a4paper,
UKenglish,
envcountsame,
]{llncs}

\usepackage[utf8]{inputenc}
\usepackage[textsize=scriptsize,obeyFinal]{todonotes}

\usepackage[T1]{fontenc}
\usepackage{color}
\usepackage{xcolor}

\usepackage{amsmath, amsfonts, amssymb, dsfont}

\usepackage{enumerate}

\usepackage{thmtools}
\usepackage{wrapfig}
\usepackage{float}

\def\ta{\mathtt{a}}
\def\tb{\mathtt{b}}
\def\tc{\mathtt{c}}

\def\tx{\mathtt{x}}

\def\ScatFact{\operatorname{ScatFact}}
\def\alph{\operatorname{alph}}
\def\r{\operatorname{r}}
\def\m{\operatorname{m}}

\def\nth#1{#1$^{\text{th}}$}

\def\abs#1{|#1|}

\renewcommand{\epsilon}{\varepsilon}
\renewcommand{\phi}{\varphi}

\def\N{\mathbb{N}} 

\DeclareMathOperator{\Pref}{Pref}
\DeclareMathOperator{\Suff}{Suff}

\DeclareMathOperator{\PerfUniv}{PUniv}
\DeclareMathOperator{\NUniv}{NUniv}
\DeclareMathOperator{\Univ}{Univ}
\DeclareMathOperator{\inner}{in}

\DeclareMathOperator{\cond}{cond}
\DeclareMathOperator{\ar}{ar}

\usepackage{bbm}

\include{macros-wpnw}
\usetikzlibrary{positioning, shapes,patterns,shadows.blur, arrows,decorations,arrows,automata,shadows,patterns,chains,graphs,calc,intersections,matrix,fit,shapes,chains,decorations.pathreplacing}
\usepackage{paralist}
\usepackage{algorithm2e}

\definecolor{cau}{RGB}{156,10,125} 

\spnewtheorem*{appendixlemma}{Lemma}{\bfseries}{}

\newif\ifpaper
\paperfalse 

\title{$m$-Nearly $k$-Universal Words - Investigating Simon Congruence}
\titlerunning{Nearly $k$-universal Words}

\author{Pamela Fleischmann\inst{1}\and
Lukas Haschke\inst{1} \and
Annika Huch\inst{1} \and
Annika Mayrock\inst{1} \and
Dirk Nowotka\inst{1}}

\institute{Kiel University, Germany, \email{$\{$fpa,lha$\}$@informatik.uni-kiel.de, 
stu216885@mail.uni-kiel.de, stu217133@mail.uni-kiel.de, dn@informatik.uni-kiel.de}}

\authorrunning{P. Fleischmann, \and L. Haschke, \land A. Huch \and A. Mayrock \land D. Nowotka}

\begin{document}

\maketitle

\begin{abstract}
Determining the index of the Simon congruence is a long outstanding open problem. Two words $u$ and 
$v$ are called Simon congruent if they have the same set of scattered factors, which are parts of 
the word in the correct order but not necessarily consecutive, e.g., $\mathtt{oath}$ is a scattered 
factor of $\mathtt{logarithm}$. Following the idea of scattered factor $k$-universality, we 
investigate $m$-nearly $k$-universality, i.e., words where $m$ scattered factors of length $k$ are 
absent, w.r.t. Simon congruence. We present a full characterisation as well as the index of the 
congruence for $m=1$. For $m\neq 1$, we show some results if in addition $w$ is $(k-1)$-universal
as well as some further insights for different $m$.

\end{abstract}

\section{Introduction}
Given a word $w$, a \emph{scattered factor} (also known as {\em (scattered) subsequence or subword}) 
of $w$ is a word, that is obtained by deleting letters from $w$ while preserving the order, i.e., 
formally $u$ of length $n\in\N_0$ is a scattered factor of $w$ (denoted by $u\in\ScatFact(w)$) if 
$w= v_1u[1]v_2u[2]...v_nu[n]v_{n+1}$ for existing (possibly empty) words $v_1,...,v_{n+1}$. For 
instance, $\mathtt{power,flow,poor,wow}$ are scattered factors of $\mathtt{powerflower}$ 
 but $\mathtt{rope}$, $\mathtt{loop}$ are not scattered 
factors since the  letters do not occur in the correct order in $w$. Therefore, scattered factors 
can be seen as a representation of $w$ with some lost data. Hence, scattered factors are not only of 
a theoretical interest, but a practical, too. When examining discrete data, e.g., protein sequences 
or incomplete or faulty transmissions of signals, scattered factors can be used as a representation 
(cf., \cite{do2021using,wang2020neural}). For instance, a faulty transmission may be 
reconstructed using scattered factors as described in 
\cite{fleischmann2021reconstructing,dress2005reconstructing,mavnuch2000characterization}. 
Scattered 
factors are also useful in sign language recognition \cite{li2010automatic}
and to get alignment-free phylogeny of whole genomes or with biological 
subwords to detect protein S-sulfenylation sites \cite{comin2012alignment,do2021using}.
Moreover, scattered factors can be found in some famous algorithmic problems like searching for 
longest (increasing) subsequences 
\cite{blumer1985smallest,bergroth2000survey,baik1999distribution}, shortest common 
supersequences \cite{maier1978complexity}, 
string-to-string correction problems \cite{wagner1974string}, most unusual time series 
subsequence \cite{keogh2007finding}, fast subsequence matching in time-series databases 
\cite{faloutsos1994fast}. Furthermore, there exist neural machine translations, which use rare 
words with subword units \cite{sennrich2015neural} or byte-level subwords \cite{wang2020neural}.

The foundations of scattered factors research were introduced by Higman \cite{higman1952ordering}, 
where it is shown that an infinite set of words
always contains words $u$ and $w$ with $u\in\ScatFact(w)$. Applications of these results can be 
found in \cite{haines1969free}. In 1972, Simon
defined the famous congruence relation regarding scattered factors in the context of piecewise testable events \cite{Simon72}, today known as {\em Simon congruence}: two words $x$ and $y$ are called {\em congruent w.r.t. $k\in\N$} ($x \sim_k y)$, iff $x$ and $y$ have the same set of scattered factors of length $k$, i.e., $\ScatFact_k(x)=\ScatFact_k(y)$ with the index denoting the length of the considered scattered factors. Thus,
we have $\mathtt{aba} \sim_2 \mathtt{aabaa}$ since $\operatorname{ScatFact}_2(\mathtt{aba}) = \{\mathtt{aa}, \mathtt{ab}, \mathtt{ba}\} = \operatorname{ScatFact}_2(\mathtt{aabaa})$ and $\mathtt{aba} \not \sim_2 \mathtt{abab}$ since $\mathtt{bb}$ is a scattered factor of $\mathtt{abab}$ but not of $\mathtt{aba}$. A profound introduction into scattered factors and Simon congruence can be found in \cite[Section 6]{lothaire} by Sakarovich and Simon.

Although $\sim_k$ is well studied from different perspectives with deep insights (cf. 
\cite{Simon72,lothaire,KufMFCS}), determining its index, i.e., determining $|\Sigma^{\ast}/\sim_k|$ 
for a given alphabet $\Sigma$ and given $k\in\N$, is 
still an open problem. First, in \cite{barker2020scattered} a special class of words, the {\em 
$k$-universal words}, were 
investigated. A word is called {\em $k$-universal} if $\ScatFact_k(w)=\Sigma^k$, i.e., $w$ has all 
the possible scattered factors of length $k$. By the definitions of $\sim_k$, we have that all 
these words are in one congruence class. These words were further investigated and characterised in 
 \cite{day2021edit,barker2020scattered,fleischmann2021scattered}. Notice that 
the idea of $k$-universality coincides with the notion of $k$-richness (cf. 
\cite{barker2020scattered} for explanations) investigated in the context of piecewise testable 
languages \cite{KarandikarKS15,CSLKarandikarS,journals/lmcs/KarandikarS19}.  One of the 
main insights of $k$-universal words is
that a word $w$ is $k$-universal iff $w$'s arch factorisation \cite{hebrard1991algorithm} has $k$ 
arches. 

Pursuing the idea of $k$-universality, where the main focus is on the cardinality of a word's  
scattered factors set rather than on the question whether two words are congruent, one can define 
the sets $M_{i,k}=\{L\subseteq\Sigma^{\ast}|\,\exists w\in\Sigma^{\ast}:\,\ScatFact_k(w)=L,|L|=i\}$ 
for
all $1\leq i\leq|\Sigma|^k$, i.e. $M_{i,k}$ contains all languages of cardinality $i$ which occur as
a scattered factor set of some word $w$ w.r.t. a length $k$. Notice that each such $L$ is a 
congruence class of $\sim_k$ and $M_{|\Sigma|^k,k}=\{\Sigma^k\}$ is built by the $k$-universal 
words. In this work, we investigate the sets $M_{i,k}$ for $i<|\Sigma|^k$. Since our main results
are for words where exactly one scattered factor from the possible scattered factor set is absent,
we call a word $m$-nearly $k$-universal if $|\ScatFact_k(w)|=|\Sigma|^k-m$, i.e. $k$-universal words
are $0$-nearly $k$-universal in the new notion.  For instance, the 
word $\mathtt{aabb}$ is $1$-nearly $2$-universal since $\tb\ta$ is absent and $\mathtt{aab}$ is
$2$-nearly $2$-universal since $\tb\ta$ and $\tb\tb$ are absent. A special subclass of $m$-nearly 
$k$-universal words has recently been studied from an algorithmic point of view in 
\cite{kosche2021absent}. There the authors investigated shortest absent scattered factors of words, 
i.e., for a given $(k-1)$-universal word $w$ the set of words with length $k$ that are not scattered 
factors of $w$. If this set has cardinality $m$, we obtain a subset of $m$-nearly $k$-universal 
words. This subset may be proper since there exists words with $m$ absent scattered factors of 
length $k$ without $w$ being $(k-1)$-universal, witness by the word $\ta\ta\tb\tb\tb$ which is 
$13$-nearly $4$-universal but not $3$-universal.

\textbf{Our contribution.} In this work, we give a full characterisation of $1$-nearly 
$k$-universal words as well as all congruence classes occurring in this subset of $\Sigma^{\ast}$.
The latter result is obtained by an algorithm that computes in linear time for a given $u$ of a 
length $k$ a word $w$ such that $u$ is the only absent scattered factor of $w$. Moreover,
we present an algorithm which decides in linear time whether a word is $1$-nearly $k$-universal.
Afterwards, we give some first insights into $m$-nearly $k$-universal words for $m>1$. Our main 
result in this part is built on the algorithm in \cite{kosche2021absent} by putting this 
algorithmic result into a combinatorial context, i.e. we are able to determine the number of absent 
scattered factors and giving the congruence classes w.r.t. $\sim_k$ for these sets.

\textbf{Structure of the work.}
In Section~\ref{prelims} we give the basic definitions and notations regarding scattered factors 
and $m$-nearly $k$-universality. In Section~\ref{nearly} we present the results on $1$-nearly 
$k$-universal words including the characterisation and the congruence classes wr.r.t. $\sim_k$. The 
result for $m>1$ are presented in Section~\ref{mnearly}. 

\section{Preliminaries}\label{prelims}
\label{sec:prelims}
Let $\mathbb{N}$ be the set of all natural numbers, $\mathbb{N}_0 = \mathbb{N} \cup \{0\}$, $[n] = \{1,\ldots,n\}$, and $[n]_0 := [n] \cup 0$. 

An \emph{alphabet} $\Sigma$ is a non empty finite set whose elements are called \emph{letters}. Set 
$\sigma=|\Sigma|$.  A \emph{word} is a finite sequence of letters from $\Sigma$. Let $\Sigma^*$ be 
the set of all finite words over $\Sigma$ with concatenation and the empty word $\epsilon$ as 
neutral element.
Set $\Sigma^+ := \Sigma^* \setminus \{\epsilon\}$. Let $w\in\Sigma^{\ast}$. 
For all $n\in\N_0$ define inductively, $w^0=\varepsilon$ and $w^n=ww^{n-1}$.
The {\em length} of $w$ is the number of $w$'s letters; thus $|\epsilon| = 0$. For all $k \in 
\mathbb{N}_0$ set $\Sigma^k := \{w \in \Sigma^* \mid |w| = k\}$ and denote $w$'s \nth{$i$} letter by $w[i]$ 
and by $w[i..j]$ denote $w[i]\cdots w[j]$ if $i < j$, $w[i]$ if $i=j$, and $\epsilon$ if $ i > j$ for all $i,j \in [\vert w \vert]$. 
Set $\alph(w) = \{\ta \in \Sigma \mid \exists i \in [\abs{w}]: w[i] = \ta \}$ as $w$'s alphabet and 
for each $\ta \in \Sigma$ set $|w|_{\ta} = |\{i \in [|w|] \mid w[i]=\ta \}|$. The word $u \in 
\Sigma^*$ is called a \emph{factor} of $w$ if there exist $x,y \in \Sigma^*$ such that $w = xuy$.
In the case $x=\epsilon$, we call $u$ a \emph{prefix} of $w$ and \emph{suffix} if $y = \epsilon$.
Let $\operatorname{Fact}(w)$, $\operatorname{Pref}(w)$ and $\operatorname{Suff}(w)$, respectively, be the sets of all factors, prefixes and suffixes of $w$.
Define the \emph{reverse} of $w$ by $w^R = w[|w|]\cdots w[1]$ and if 
$w=x_1^{k_1}x_2^{k_2}\cdots x_{\ell}^{k_{\ell}}\in\Sigma^{\ast}$ with 
	$k_i,\ell\in\N$, $i\in[\ell]$, the {\em condensed form (print)} of $w$ is defined by $\cond(w)=x_1\cdots 
	x_{\ell}$ assumed that $x_j\neq x_{j+1}$ for $j\in[\ell-1]$. 
Let $<_{\Sigma}$ be a total order on $\Sigma$. We extend this order to the \emph{lexicographical 
order} on $\Sigma^*$ by $u < v$ for $u,v \in \Sigma^*$ iff $u \in \operatorname{Pref}(v)$ or $u = 
x\ta u'$ and $v = x \tb v'$ with $\ta < \tb$ for $\ta , \tb \in \Sigma$ and some $u',v',x \in 
\Sigma^*$. Define $w_{\Sigma}$ as the word in $\Sigma^{\sigma}$ with 
$w_{\Sigma}[i]<_{\Sigma}w_{\Sigma}[i+1]$ and $\alph(w)=\Sigma$.
For further definitions see \cite{lothaire}.

After fixing the basic notations, we introduce the scattered factors.

\begin{definition}
	Let $w \in \Sigma^*$ and $n \in \mathbb{N}_0$. A word $u \in \Sigma^n$ is called a \emph{scattered factor} of w ($u \in \ScatFact(w)$) if there exist $v_1,\dots, v_{n+1} \in \Sigma^*$ such that $w = v_1 u[1] v_2 u[2] \cdots v_n u[n]$ $v_{n+1}$.
	Set $\ScatFact_k(w)=\{u\in\ScatFact(w)|\,|u|=k\}$.
\end{definition}

To give an example $\mathtt{cau}, \mathtt{flower}, \mathtt{cafe}$, $\mathtt{life}$ and $\mathtt{ufo}$ are all scattered factors of $\mathtt{cauliflower}$ but neither $\mathtt{flour}$ nor $\mathtt{row}$.

Tightly related to the notion of scattered factors is the famous {\em Simon congruence}. Two words 
are congruent modulo $k\in\N_0$ if they have the same set of scattered factors of length $k$, e.g., 
$\mathtt{aba}$ and $\mathtt{aabaa}$ are congruent w.r.t. $2$ since 
$\operatorname{ScatFact}_2(\mathtt{aba}) = \{\ta\ta, \ta\tb, \tb\ta\} = 
\operatorname{ScatFact}_2(\mathtt{aabaa})$. 

\begin{definition}
	Two words $w,v \in \Sigma^*$ are \emph{Simon congruent} w.r.t. $k\in\N_0$ ($w \sim_k v$) if $\operatorname{ScatFact}_k(w) = \operatorname{ScatFact}_k(v)$.
\end{definition}

Since $\ScatFact_k(w)\subseteq\Sigma^k$ holds for all $k\in\N_0$, determining the index of the Simon 
congruence can be split into the parametrised problem on determining how many scattered factor sets 
- or equivalently how many different words - exist with $|\ScatFact_k(w)|=\sigma^k-m$ for all 
$m\in\N_0$. In \cite{barker2020scattered,fleischmann2021scattered,day2021edit}
the scattered factor universality was investigated, which describes the problem for $m=0$.

\begin{definition}
	A word $w \in \Sigma^*$ is called \emph{$k$-universal} if $\vert \ScatFact_k(w) \vert = 
\sigma^k$.
Let $\iota(w)$ denote the {\em universality index}, i.e. the largest $k\in\N_0$ such that $w$ is 
$k$-universal.
	We call a 1-universal word $w$ just \emph{universal}. Denote by 
$\operatorname{Univ}_{\Sigma,k}$ the set of all words $w$ with $\iota(w)=k$. 
	
\end{definition}

\begin{remark}{}
	By definition, all $k$-universal words are congruent modulo $k$ and a $k$-universal word $w \in \Sigma^*$ is also $k'$-universal for all $k' \leq k$.
\end{remark}

In this work, we are investigating $m$-nearly $k$-universal words. These are words, where in 
comparison to $\Sigma^k$, $m$ words of length $k$ are absent from the scattered factor set. A 
special case of these words was investigated in 
\cite{kosche2021absent} where the {\em shortest absent scattered 
factors} of a word are determined. 	In the unary alphabet $\epsilon$ is the only word which has 
$\vert \Sigma \vert^k -1 = 0$ scattered factors and the notion is not well-defined for $m>1$. 
Therefore we only consider at least binary alphabets. Moreover, we assume $\Sigma=\alph(w)$ 
for a given $w$, if not stated otherwise.
	
\begin{definition}
Let $k,m\in\N_0$.
	A word $w \in \Sigma^*$ is called \emph{$m$-nearly $k$-universal} if $\vert \ScatFact_k(w) \vert 
= \sigma^k -m$. Let $\NUniv_{\Sigma,m,k}$ denote the set of all $m$-nearly $k$-universal words in 
$\Sigma^*$. We call a $1$-nearly $k$-universal word simply {\em nearly $k$-universal}.
\end{definition}

\begin{remark}{}
	Unlike the $k$-universality, $w\in\NUniv_{\Sigma,m,k}$ does not imply 
$w\in\NUniv_{\Sigma,m,k-1}$ for $m>0$: we have $\ta\tb\ta\in\NUniv_{\Sigma,1,2}$ but by 
$\iota(\ta\tb\ta)=1$, $\ta\tb\ta$ is not $m$-nearly $1$-universal for all $m>0$. 
\end{remark}

One of the main tools for the investigation of $m$-nearly $k$-universal words is the \emph{arch factorisation} which was introduced by Hebrard \cite{hebrard1991algorithm}. In this factorisation a word is factorised into universal factors and a rest.
	
\begin{definition}
	For a word $w \in \Sigma^*$ the \emph{arch factorisation} is given by $w = \ar_1(w) \cdots \ar_k(w) \r(w)$ for $k \in \mathbb{N}_0$ with\\
	(a) $\iota(\ar_i(w))=1$ for all $i \in [k]$, \\
	(b) $\ar_i(w)[\vert \ar_i(w) \vert] \notin \alph(\ar_i(w)[1 \cdots \vert \ar_i(w) \vert - 1 ])$ for all $i \in [k]$, and\\
	(c) $\alph(\r(w)) \subset \Sigma$.\\
	The words $\ar_i(w)$ are called \emph{arches} of $w$ and $\r(w)$ is the \emph{rest} of $w$.
	Define the \emph{modus} $\m(w) = \ar_1(w)[\vert\ar_1(w)\vert] \cdots\ar_k(w)[\vert\ar_k(w)\vert]$.
	The \emph{inner of the $i^{th}$ arch of $w$} is defined as the prefix of 
$\ar_i(w)$ such that $\ar_i(w) = \inner_i(w) 
\m(w)[i]$ holds.
\end{definition}

To visualise the arch factorisation in explicit examples we use brackets. For example we write $(\mathtt{aab})\cdot (\mathtt{bba}) \cdot \ta$ to mark the two arches, namely $\mathtt{aab}$ and $\mathtt{bba}$ and the rest, $\ta$, which is denoted without brackets.
	
\begin{remark}{}
	The modus $\m(w)$ consists of all unique last letters of the arches and is therefore uniquely defined.
\end{remark}

Based on the arch factorisation we define perfect universal words, which are words without a rest.

\begin{definition}
	We call a word $w \in \Sigma^*$ \emph{perfect $k$-universal} if $\iota(w)=k$ and $\r(w) = \epsilon$. The set of all these words with $\alph(w)=\Sigma$ is denoted by $\PerfUniv_{\Sigma,k}$.
\end{definition}

For the algorithmic results in Section~\ref{nearly} and \ref{mnearly} we use the standard 
computational model RAM with logarithmic word-size (see, e.g., \cite{KarkkainenSB06}), i.e., we 
follow a standard assumption from stringology, if $w$ is the 
input word for our algorithms, we assume $\Sigma=\alph(w)=\{ 1,2, \ldots, \sigma\}$.

\section{Nearly $k$-Universal Words}\label{nearly}
\label{sec:arbitrary}
In this section we characterise the nearly $k$-universal words, i.e. words $w\in\Sigma^{\ast}$ with $\ScatFact_k(w)=\sigma^k-1$ for a fixed $k\in\N$. 
Moreover, we show that there exist exactly $\sigma^k$ different classes w.r.t. $\sim_k$, i.e., for 
each word $v\in\Sigma^k$ there exists a word $w\in\Sigma^{\ast}$ such that 
$\ScatFact_k(w)=\Sigma^k\backslash\{v\}$. First, we show that some peculiarities do not 
occur for $m=1$. Consider the words $w = \mathtt{(accab)}, w' = \mathtt{(abc) \cdot a} \in 
\NUniv_{\Sigma,3,2}$. Notice that $r(w)=\varepsilon$ and $|\alph(r(w'))|<\sigma-1$. These cases
cannot occur for $m=1$. Also, in general, $\iota(w)<k-1$ is possible,
witnessed by $\mathtt{(abac) \cdot bac} \in \NUniv_{\Sigma,7,3}$, but it is not possible for nearly 
$k$-universal words. 

\begin{theorem}{}\label{theorem_univ}
	If $w\in\NUniv_{\Sigma,1,k}$  then $\iota(w) = k-1$ and $\vert\alph (\operatorname{r}(w)) \vert = \sigma - 1$.
\end{theorem}
\ifpaper
\else    
\begin{proof}
Suppose $\iota(w) < k-1$. Choose $v \in \Sigma^{k-1}$ with $v \notin \ScatFact_{k-1}(w)$. Thus, for all $x \in \Sigma$, $vx \notin \ScatFact_{k}(w)$ and we obtain 
$\vert \ScatFact_{k}(w) \vert \leq \vert \Sigma \vert^k - \vert \Sigma \vert < \vert \Sigma \vert^k - 1$ - a contradiction. By $\iota(w)=k-1$, we have $\vert \operatorname{m}(w) \vert = k-1$. If $w\in\PerfUniv_{\Sigma,k-1}$, neither $\m(w)\ta$ nor $\m(w)\tb$ were scattered factors of $w$, for $\ta\neq\tb$, and $w$ would not be nearly $k$-universal. The same argumentation holds if there exists $\ta,\tb\in\Sigma\backslash\alph(\r(w))$ different. This concludes the proof.\qed
\end{proof}

\fi

\begin{remark}{} \label{alph}
	Theorem~\ref{theorem_univ} implies 	that the length a of nearly $k$-universal word is at least  
$k\sigma - 1$ since we have 
	$k-1$ arches and a rest of length $\sigma-1$. Moreover, for each nearly $k$-universal word $w$ exists a unique letter $\ta_w$ with 
	$\alph(\r(w))=\Sigma\backslash\{\ta_w\}$.
\end{remark}

\begin{corollary}{}\label{absent}
    Given $w\in\NUniv_{\Sigma,1,k}$, we have $\ScatFact_k(w)=\Sigma^k\backslash\{\m(w)\ta_w\}$.
\end{corollary}
\ifpaper
\else    
\begin{proof}
    Follows directly by Theorem \ref{theorem_univ}.\qed
\end{proof}

\fi

The conditions of Theorem~\ref{theorem_univ} do not suffice for a characterisation of nearly 
$k$-universal words. Consider the word $w = (\mathtt{acb}) \cdot \tb\ta$ with 
$\iota(w)=1$ and $\alph(\r(w)) = \{\ta,\tb\}$. We have 
$|\ScatFact_2(w) \vert = |\Sigma^2\backslash\{\tc\tc,\tb\tc\}|$ and thus 
$w\not\in\NUniv_{\Sigma,1,2}$. The first, na\"ive characterisation
uses Corollary~\ref{absent}: all words of length $k$ ending in $\ta_w$, but $\m(w)\ta_w$, have to 
appear within the word (all others appear necessarily).

\begin{proposition}\label{pidgeon}
	A word $w \in \NUniv_{\Sigma,1,k}$ iff $\iota(w) = k-1$, $\alph(\r(w))=\Sigma\backslash\{\ta_w\}$, and for all 
	$v\in\Sigma^k$ with $v[1..k-1]\neq \m(w)$ and $v[k]=\ta_w$ there exists $i \in [k-2]$ with $v[i]v[i+1] \in \operatorname{ScatFact}_2(\operatorname{ar}_i(w))$ or
			$v[k-1]\ta_w \in \operatorname{ScatFact}_2(\operatorname{ar}_{k-1}(w))$.
\end{proposition}
\ifpaper
\else    
\begin{proof}
	First, let $w$ be nearly $k$-universal. The first two claims follow immediately from Theorem \ref{theorem_univ} and Corollary~\ref{absent}. Moreover, we know 
	$\operatorname{m}(w)\ta_w \notin \operatorname{ScatFact}_k(w)$. Let $v\in\Sigma^k$ with $v[1..k-1]\neq \m(w)$ and $v[k]=\ta_w$.
	
	If $v[k-1]\ta_w\in\ScatFact_2(\ar_{k-1}(w))$, we are done. Thus, assume that $v[k-1]\ta_w\not\in
	\ScatFact_2(\ar_{k-1}(w))$. Since $w\in\NUniv_{\Sigma,1,k}$, $v\neq \m(w)\ta_w$, and $\ta_w\not\in\alph(r(w))$, we get immediately by the pidgeon hole principle that there exists $i\in[k-2]$ with $v[i]v[i+1]\in\ScatFact_2(\ar_i(w))$.

	Second assume the three constraints to hold true and suppose that $w$ is not nearly $k$-universal. We get immediately $\m(w)\ta_w \notin \operatorname{ScatFact}_k(w)$. Since by supposition $w$ is not nearly $k$-universal there exists $u\not\in\ScatFact_k(w)$ with $u\neq \m(w)\ta_w$.
	Since $\iota(w)=k-1$, we have for all $u\in\Sigma^{k-1}$ and $\ta\in\alph(r(w))$ immediately $u\ta\in\ScatFact_k(w)$. This implies $u[k]=\ta_w$ and $u[1..k-1]\neq m(w)$. Thus, there exists $i\in[k-2]$ with $u[i]u[i+1]\in\ScatFact_2(\ar_i(w))$ or $u[k-1]\ta_w\in\ScatFact_2(\ar_{k-1}(w))$. In the first case
	$u[1..i-1]$ can be chosen from the first $i-1$ arches (which build a $(i-1)$-universal word) and $u[i+2..k-1]\ta$ can be chosen from $\ar_{i+1}(w)\cdots
	\ar_{k-1}(w)$ (which is a $(k-i-1)$-universal word and $|u[i+2..k-1]|=k-i-1$). In the second case $u[1..k-2]$ be chosen from the first $k-2$ arches.
	Thus, in both cases we have $u\in\ScatFact_k(w)$ - a contradiction.\qed

\end{proof}

\fi

Since $\ta\tc,\tc\tc\in\ScatFact_2(\ar_1(w))$ and $\ta\tc,\tb\tc\in\ScatFact_2(\ar_2(w))$,
we have $w=\mathtt{(accb)\cdot (bac) \cdot ab} \in \NUniv_{\Sigma,1,3}$.
This characterisation is not very helpful since checking whether a word is nearly $k$-universal 
means to check all $\sigma^{k-1}$ options for $v$.
The following characterisation does not only provide an efficient way to check whether 
$w\in\NUniv_{\Sigma,1,k}$ but also builds the basis for an
efficient algorithm regarding $\sim_k$. In beforehand, we 
prove that {\em cutting off} $\ell$ arches at the beginning of a nearly $k$-universal word, leads to 
a nearly $(k-\ell)$-universal word.

\begin{lemma}\label{cutting_first}
	Let $\ell\leq k-1$. If $w\in\NUniv_{\Sigma,1,k}$ with  $w = \ar_1(w)\cdots \ar_{k-1}(w)\r(w)$, then $\ar_{\ell+1}(w) ... \ar_{k-1}(w)\r(w)\in\NUniv_{\Sigma,1,k-\ell}$.
\end{lemma}
\ifpaper
\else
\begin{proof}
   It suffices to prove the claim for $\ell=1$; the main statement follows inductively. Set $w'=\ar_2(w)\cdots\ar_{k-1}(w)$. Let $\mathtt{x} = \operatorname{m}(w)[1]$. We have exactly $\frac{1}{\sigma} \cdot \sigma^k$ scattered factors of length $k$ with first letter $\mathtt{a}$ for all $\ta \in \Sigma \setminus \{\mathtt{x}\}$. The number of scattered factors starting with $\mathtt{x}$ is 
   $\frac{1}{\sigma} \cdot \sigma^k - 1$ because $\operatorname{m}(w)\mathtt{a}_w$ is not a scattered factor of $w$. 
   Set $S_\ta=\{u|\,\ta u\in\ScatFact(w)\}$ for all $\ta\in\Sigma$. Thus, for all $\ta,\tb\in\Sigma\backslash\{\tx\}$ we have
   $S_\ta=S_\tb$. Moreover, we have $S_\tx=S_\ta\backslash\{\m(w)[2..k-1]\ta_w\}$. By $x=\m(w)[1]$ follows that $u\in\ScatFact_{k-1}(w')$ for all $u\in S_\tx$. These are $\sigma^{k-1}-1$ many and by $\m(w)[2..k-1]\ta_w\not\in\ScatFact_{k-1}(w')$, the claim is proven.\qed
\end{proof}

\fi

\begin{remark}{}
 Notice that Lemma~\ref{cutting_first} is not applicable for arches in the middle:  
 $(\mathtt{ab}) \cdot (\mathtt{aa}\mathtt{b}) \cdot 
\mathtt{b}\in\NUniv_{\Sigma,1,3}$ but $(\ta\tb)\cdot \tb\not\in\NUniv_{\Sigma,1,2}$. 
Moreover, Lemma~\ref{cutting_first} does not hold for $m>1$: $w = \mathtt{(abc) \cdot 
(bca) \cdot bb} \in \NUniv_{\Sigma,7,3}$ but $ 
\mathtt{(bca) \cdot bbb} \in \NUniv_{\Sigma,3,2} \neq 
\NUniv_{\Sigma,7,2}$.
\end{remark}

Now we present a more suitable charactersiation for nearly $k$-universal words.
Here, $\ScatFact_k(w^R)=\{u^R|\,u\in\ScatFact_k(w)\}$ plays an important role.

\begin{theorem}{}\label{niceequivalence}
For $w\in\Sigma^{\ast}$ the following statements are equivalent\\
(1) $w \in \NUniv_{\Sigma,1,k}$,\\
(2) $\iota(w) =  k-1$, $\vert \alph(\r(w))\vert = \sigma -1=\vert \alph(\r(w^R))\vert$, 
and \\
(a) if $k$ is even then there exists $u_1,v_2 \in \PerfUniv_{\Sigma,\frac k 2}$,  $u_2,v_1 \in 
\PerfUniv_{\Sigma,\frac k 2 -1}$ and $x_i \in \Sigma^+$ with $\vert\alph (x_i)\vert = \sigma 
-1$ with $w = u_ix_i v_i^R$ for $i \in [2]$.\\
(b) if $k$ is odd then there exist $u,v \in \PerfUniv_{\Sigma,\frac {k-1} 2}$, and $x \in \Sigma^+$ 
with $\vert\alph (x)\vert = \sigma -1$ with $w = uxv^R$.\\
(3)  $\iota(w) =  k-1$, $\vert \alph(\r(w))\vert = \sigma -1=\vert \alph(\r(w^R))\vert$, and
	for all $\hat{k},\tilde{k} \in \mathbb{N}$ with $\hat{k} + \tilde{k} +1 = k$ there exist $u \in 
\PerfUniv_{\hat{k}}$, $v \in \PerfUniv_{\tilde{k}}$, and $x \in 
\Sigma^+$ with $\vert \operatorname{alph}(x) \vert = \sigma -1$ such that $w = uxv^R$.
\end{theorem}
\ifpaper
\else    
\begin{proof}
First, we prove (2) implies (1).    
	We have to show that $w$ is nearly $k$-universal under the three constraints.
	We know $\m(w)\ta_w\not\in\ScatFact_k(w)$. Let $y\in\Sigma^k\backslash\{\m(w)\ta_w\}$. If $y[k]\neq\ta_w$,
	we have immediately $y\in\ScatFact_k(w)$ by the second condition. Thus, assume $y[k]=\ta_w$.\\
	\textbf{case 1:} $k$ is even \\
	Choose $u_1,u_2,v_1,v_2,x_1,x_2$ according to condition $a)$. Since $u_1$ and $v_2$ are perfect $\frac{k}{2}$-universal and $u_2$ and $v_1$ are perfect $\frac{k}{2}-1$-universal, we have 
	\begin{align*}
	y\left[1..\frac{k}{2}\right] &\in\ScatFact(u_1), & y\left[\frac{k}{2}+2..k\right]\in\ScatFact(v_1^R),\\
	y\left[1..\frac{k}{2}-1\right]  &\in\ScatFact(u_2), & y\left[\frac{k}{2}+1..k\right]\in\ScatFact(v_2^R).
	\end{align*}
	Thus, if $y[\frac{k}{2}+1]\in\alph(x_1)$ or $y[\frac{k}{2}]\in\alph(x_2)$, we have $y\in\ScatFact_k(w)$.
	Assume $y[\frac{k}{2}+1]\not\in\alph(x_1)$ and $y[\frac{k}{2}]\not\in\alph(x_2)$. Since we have also proven the claim if two consecutive letters of $y$ are in one arch of $u_1$, $u_2$, $v_1$, or $v_2$, we may assume
	that $y[1..\frac{k}{2}-1]=\m(u_2)$ and $y=[\frac{k}{2}+2..k]=\m(v_1)^R$. By $|\alph(x_1)|=\sigma-1$, we have
	$y[\frac{k}{2}+1]=\m(v_2)[\frac{k}{2}]$, and analogously by $|\alph(x_2)|=\sigma-1$, we have $y[\frac{k}{2}]=
	\m(u_1)[\frac{k}{2}]$. Choose $i_1,i_2\in[|w|]$ with $w[i_1]=\m(v_2)[\frac{k}{2}]$ and $w[i_2]=\m(u_1)[\frac{k}{2}]$.
	If $i_1\geq i_2$, $w$ would have at least $\iota(u_1)+\iota(v_2^R)=k$ arches - a contradiction. Thus we have
	$i_1<i_2$. This implies that $y[\frac{k}{2}+1]$ has be chosen before $y[\frac{k}{2}]$ in $w$. This implies $\iota(w)<k-1$ - a contradiction\\
	
	\begin{figure}
	\label{niceequivalence_example}
\begin{tikzpicture}[scale=0.5,
blackbox/.style={rectangle, draw=black, thick, minimum height=0.6cm},
redbox/.style={rectangle, draw=red!75, 
			 thick, minimum height=0.5cm}
]

\coordinate (u1) at (10,1);
\coordinate (v1) at (12,1);
\coordinate (u2) at (8.5,1);
\coordinate (v2) at (10.5,1);
\coordinate (e) at (20,1);

\draw[-] (0,0.75) -- (0,1.25);
\draw[-] (20,0.75) -- (20,1.25);

\node (w) at (-1,1) {$w$};
\draw[-] (u1) -- ($ (u1)+(0,0.25) $);
\draw[-] (v1) -- ($ (v1)+(0,0.25) $);

\draw[-] (u2) -- ($ (u2)+(0,-0.25) $);
\draw[-] (v2) -- ($ (v2)+(0,-0.25) $);

\draw[-] (0,1) --  node[above] {$u_1$} (u1) 
		-- node[above] {$x_1$} (v1)
		-- node[above] {$v_1^R$} (e);
		
\draw[-] (0,1) --  node[below] {$u_2$} (u2) 
		-- node[below] {$x_2$} (v2)
		-- node[below] {$v_2^R$} (e);
		
\draw[gray] (0.1,1.7) to[out=90,in=90, distance=1cm] node[above, gray, font=\scriptsize] {$\ar_1$} (1.5,1.75);
\draw[gray] (1.5,1.75) to[out=90,in=90, distance=1cm] node[above, gray, font=\scriptsize] {$\ar_2$} (3,1.75);
\draw[gray] ($ (u1)+(-1.5,0.75) $) node[left] {$\ldots$} to[out=90,in=90, distance=1cm] node[above, gray, font=\scriptsize] {$\ar_\frac{k}{2}$} ($ (u1)+(-0.1,0.75) $);
\draw[gray] ($ (v1)+(0.1,0.75) $) to[out=90,in=90, distance=1cm] node[above, gray, font=\scriptsize] {$\ar_1$} ($ (v1)+(1.5,0.75) $);
\draw[gray] ($ (e)+(-1.5,0.75) $) node[left] {$\ldots$} to[out=90,in=90, distance=1cm] node[above, gray, font=\scriptsize] {$\ar_{\frac{k}{2}-1}$} ($ (e)+(-0.1,0.75) $);
\draw[gray] (0.1,0.25) to[out=-90,in=-90, distance=1cm] node[below, gray, font=\scriptsize] {$\ar_1$} (1.5,0.25);
\draw[gray] ($ (u2)+(-1.5,-0.75) $) node[left] {$\ldots$} to[out=-90,in=-90, distance=1cm] node[below, gray, font=\scriptsize] {$\ar_{\frac{k}{2}-1}$} ($ (u2)+(-0.1,-0.75) $);
\draw[gray] ($ (v2)+(0.1,-0.75) $) to[out=-90,in=-90, distance=1cm] node[below, gray, font=\scriptsize] {$\ar_1$} ($ (v2)+(1.5,-0.75) $);
\draw[gray] ($ (v2)+(1.6,-0.75) $) to[out=-90,in=-90, distance=1cm] node[below, gray, font=\scriptsize] {$\ar_2$} ($ (v2)+(3,-0.75) $);
\draw[gray] ($ (e)+(-1.5,-0.75) $) node[left] {$\ldots$} to[out=-90,in=-90, distance=1cm] node[below, gray, font=\scriptsize] {$\ar_\frac{k}{2}$} ($ (e)+(-0.1,-0.75) $);

\node[below, gray, font=\tiny] (y1) at ($ (1.5,1)+(-0.35,0) $) {$y[1]$};
\node[below, gray, font=\tiny] (yk2-1) at ($ (u2)+(-0.35,0) $) {$y[\frac{k}{2}-1]$};
\node[above, gray, font=\tiny] (yk2) at ($ (u1)+(-0.35,0) $) {$y[\frac{k}{2}]$};
\node[below, gray, font=\tiny] (yk2+1) at ($ (v2)+(0.5,0) $) {$y[\frac{k}{2}+1]$};
\node[above, gray, font=\tiny] (yk2+2) at ($ (v1)+(0.5,0) $) {$y[\frac{k}{2}+2]$};
\node[above, gray, font=\tiny] (yk) at ($ (e)+(-1.25,0) $) {$y[k]$};

\end{tikzpicture}
	\caption{The factorisation of $w$ for even $k$ where $y$'s letters occur as the 
modus.}
	\end{figure}
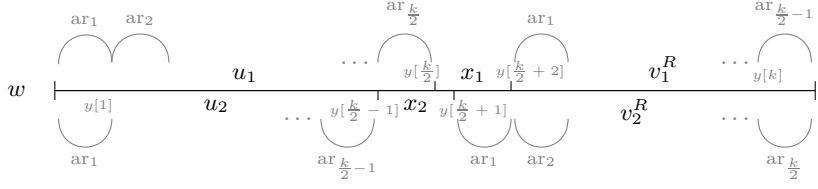

	\textbf{case 2:} $k$ is odd\\
	Choose $u,v,x$ according to condition $b)$. Since $u$ and $v$ are perfect $\frac{k-1}{2}$-universal, we have 
	\begin{align*}
	y\left[1..\frac{k-1}{2}\right] &\in\ScatFact(u), & y\left[\frac{k-1}{2}+2..k\right]\in\ScatFact(v).
	\end{align*}
	If $y[\frac{k-1}{2}+1]\in \alph(x)$, the claim is proven. Thus, assume $y[\frac{k-1}{2}+1]\not\in\alph(x)$. Again
	we can also assume $y[1..\frac{k-1}{2}]=\m(u)$ and $y[\frac{k-1}{2}+2..k]=\m(v)^R$. Since $|\alph(x)|=\sigma-1$,
	we have $y[\frac{k-1}{2}+1]=\m(w)[\frac{k-1}{2}+1]$ which occurs after $\m(v)[\frac{k-1}{2}]$ in $w$. Again we obtain $\iota(w)<k-1$, a contradiction.
	
	\bigskip
	
Now we prove that (1) implies (3).
Consider firstly $w$ be nearly $k$-universal. Then the first two claims follow immediately by Theorem \ref{theorem_univ} and the fact that $w^R$ is nearly $k$-universal.\\
	By $\iota(w)=k-1$, we have $w,w^R\in\Univ_{\Sigma,k'}$ for all $k'\leq k-1$. Let $\hat{k},\tilde{k} \in \N_{<k}$ with $\hat{k} + \tilde{k} + 1 = k$. Thus, there exist $u \in \PerfUniv_{\Sigma,\hat{k}}$ and $v \in \PerfUniv_{\Sigma,\tilde{k}}$ with $u \in \operatorname{Pref}(w)$ and $v \in \operatorname{Pref}(w^R)$.
	Choose $x\in\Sigma^{\ast}$ with $w = u x v^R$.
	By Lemma~\ref{cutting_first}, we get 
	$x v^R\in\NUniv_{\Sigma,1,k-\hat{k}}$.
	Thus, $v x^R\in\NUniv_{\Sigma,1,k-\hat{k}}$.
	Applying Lemma~\ref{cutting_first} again, we obtain $x^R \in\NUniv_{\Sigma,1,1}$. 
	By Theorem~\ref{theorem_univ} we get $|\alph(x)|=\sigma-1$.
	
\bigskip

Since (3) implies (2) immediately, the claim is proven.\qed
\end{proof}

\fi

We have $w = (\ta\ta\tb)\cdot(\tb\ta)\cdot(\ta\tb)\cdot \ta\not\in\NUniv_{\Sigma,1,4}$ 
since we have the factorisation $(\ta\ta\tb)(\tb\ta)\cdot \ta \cdot ((\ta\tb))^R$ meeting the 
requirements but also the factorisation $(\ta\ta\tb)\cdot\varepsilon \cdot 
((\ta\tb)(\ta\ta\tb))^R$ not meeting them, witnessing that both factorisations are 
needed.

\begin{corollary}\label{palindromes}
We have $ww^R \in \NUniv_{\Sigma,1,2k-1}$ iff $w\in\NUniv_{\Sigma,1,k}$ as well as $w\ta w^R \in 
\NUniv_{\Sigma,1,2k-1}$ with $\ta \in \Sigma$ iff $w\in\NUniv_{\Sigma,1,k}$ and $\ta \in 
\alph(\r(w))$.
\end{corollary}

\ifpaper
\else
\begin{proof}
Consider first $ww^R$ with $w\in\NUniv_{\Sigma,1,k}$. Since $2k-1$ is odd, 
$|\alph(x)|=|\alph(\r(w)\r(w^R))|=\sigma-1$ 
shows that $ww^R \in \NUniv_{\Sigma,1,2k-1}$. Now let $ww^R \in \NUniv_{\Sigma,1,2k-1}$. Since 
$\vert w \vert = \vert w^R \vert$ we receive a factorisation $ww^R = uxv^R$ with $u,v \in 
\PerfUniv_{\Sigma, k-1}$, $x \in \Sigma^+$ with $\vert\alph(x)\vert = \sigma-1$ as well as $u 
\in \Pref(w)$, $v^R \in \Suff(w^R)$. Since $ww^R$ is a palindrome, we have $u = v$. Thus, $ww^R = 
uxu^R$. Moreover, $x=yy$ for  $y \in \Sigma^+$ with $\vert\alph(y)\vert = \sigma -1$ 
holds. Applying Lemma~\ref{cutting_first} with $\ell = k-1$ we obtain that $xu^R \in 
\NUniv_{\Sigma,1,k-1}$. Thus, $ux \in \NUniv_{\Sigma,1,k-1}$. Since multiple occurrences of letters 
in the rest do not have an impact on the absent scattered factor of nearly $k$-universal words, $w = 
uy, w^R = yu \in \NUniv_{\Sigma,1,k}$ follows.

Considering palindromes of odd length, we can also apply Theorem~\ref{niceequivalence} 
directly, choosing $\hat{k} = \tilde{k} = k$, we get 
$|\alph(x)|=|\alph(\r(w)\r(w^R))|=\sigma-1$. The same argumentation as in the even case proves the 
claim.  \qed
\end{proof}

\fi

With Theorem~\ref{niceequivalence} we are able to solve the following two problems (for a given 
$k$) efficiently: decide whether a word is nearly $k$-universal  and find for a given 
$u\in\Sigma^{k}$
a $w\in\NUniv_{\Sigma,1,k}$ such that $u\not\in\ScatFact_k(w)$. The latter one leads immediately to 
the index of the Simon congruence restricted to nearly $k$-universal words. Notice that for the 
first problem, a linear time algorithm is implicitly given in \cite{kosche2021absent}: if $w$ is a 
word of length $n$, the SAS tree can be constructed in time $\mathcal{O}(n)$ and in time 
$\mathcal{O}(k)$ the lexicographically smallest shortest absent scattered factors can be 
determined; if there is only one shortest absent scattered factor, we have 
$w\in\NUniv_{\Sigma,1,k}$. The following algorithm can only check whether
a word is nearly $k$-universal but therefore does not need any additional data structures.

\begin{proposition}\label{checknearly}
Given $w\in\Sigma^{\ast}$ and $k\in\N$, we can decide  whether 
$w\in\Univ_{\Sigma,1,k}$ in time $\mathcal{O}(|w|)$.
In the positive, the absent scattered factor is also computed (see Algorithm~\ref{alg:uxvR}).
\end{proposition}

\ifpaper
\else
\begin{proof}
By \cite{barker2020scattered} we know that the arch factorisation can be computed in time 
$\mathcal{O}(|w|)$.
While computing the arch factorisation of $w$, store the end of the \nth{($\frac{k}{2}-1$)} arch in $i_1$ and
the end of the \nth{$\frac{k}{2}$} in $i_2$. Analogously, while computing the arch factorisation of $w^R$, 
store the end of the \nth{($\frac{k}{2}-1$)} arch in $j_2$ and
the end of the \nth{$\frac{k}{2}$} in $j_1$. These four values can be obtained in $\mathcal{O}(|w|)$. 
Now we have to check that $w[i_1..j_2]$ and $w[i_2..j_1]$ both contain each all letters from $\Sigma$ but one.
This can be done in $\mathcal{O}(\sigma)$. By $\sigma\leq n$, the claim is proven. For $k$ odd, we only need to check
one factorisation.

While checking the conditions of Theorem~\ref{niceequivalence}, we also computed $m(w)\ta_w$ and 
thus the algorithm also determines the absent scattered factor.
\qed
\end{proof}

\fi

\ifpaper
\else
\SetKwComment{Comment}{/* }{ */}
\begin{algorithm}
	\caption{Testing nearly $k$-universality (cf. Proposition~\ref{checknearly})}\label{alg:uxvR}
	\KwData{Given $w \in \Sigma^*$ with arch factorisation and $k \in \N$.}
	\KwResult{True, if $w \in \NUniv_{\Sigma,1,k}$. False, otherwise.}
	\eIf{$\iota(w) \neq k-1 \mid\mid \vert \alph(\r(w))\vert \neq \sigma -1 \mid\mid \vert \alph(\r(w^R))\vert \neq \sigma -1$}
	{\Return{false}\;}
	{	
		\eIf{$k \mod_2 == 0$}
		{$w_{v_1} := (\ar_{\frac k 2}(w^R) \cdots \ar_{k-1}(w^R) \r(w^R))^R$ \Comment*[r]{The index denotes the deleted archs of $w$'s factorisation}
		$w_{v_2} := (\ar_{\frac k 2 +1}(w^R) \cdots \ar_{k-1}(w^R) \r(w^R))^R$\;
		\Return{$\vert \alph(\r(w_{v_1})) \vert == \sigma -1 \quad \&\& \quad \vert \alph(\r(w_{v_2})) \vert == \sigma -1$}\;
		}
		{$w_v := (\ar_{\frac {k-1} 2}(w^R) \cdots \ar_{k-1}(w^R) \r(w^R))^R$\;
		\Return{$\vert \alph(\r(w_{v})) \vert == \sigma -1$};\
		}
	}
\end{algorithm}

\fi

\begin{remark}
Theorem~\ref{niceequivalence} can also be used to construct nearly $k$-universal words: if $k$ is 
odd choose $u,v\in\PerfUniv_{\Sigma,\frac{k-1}{2}}$ as well as an $x$ with $|\alph(x)|=\sigma-1$ 
and $uxv^R$ is nearly $k$-universal. In the case that $k$ is even, choose 
$u,v\in\PerfUniv_{\Sigma,\frac{k-1}{2}}$ as well $x_1,x_2$ such that 
$|\alph(x_1)|=|\alph(x_2)|=\sigma-1$. Now, we have $ux_2yx_1v\in\NUniv_{\Sigma,1,k}$ iff 
$y[|y|]\not\in\alph(x_2)$ and $y[1]\not\in\alph(x_1)$.
\end{remark}

Now, we present an algorithm for the second problem. Please recall that 
$\Sigma_\ta=\Sigma\backslash\{\ta\}$ and $w_{\Sigma_\ta}$ is the word containing all letters of 
$\Sigma_{\ta}$ w.r.t. a predefined order $<_{\Sigma}$ on $\Sigma$.
These words can be preprocessed in time $\mathcal{O}(\sigma)$ for all $\ta\in\Sigma$.

\begin{theorem}\label{utow}
Given $u\in\Sigma^k$ for $k\in\N$, one can compute $w\in\Sigma^{\ast}$ with $\ScatFact_k(w)=\Sigma^k\backslash\{u\}$ in time $\mathcal{O}(k)$. More precisely,
there exists an algorithm needing $k$ steps computing $w\in\NUniv_{\Sigma,1,k}$ of minimal 
length (see Algorithm~\ref{alg:one}).
\end{theorem}

\ifpaper
\else
\begin{proof}
Given $k\in\N$ and $u\in\Sigma^k$ the following inductive algorithm constructs a word $w$ with $\ScatFact(w) = \Sigma^k \setminus \{u\}$:
for all  $i=[k-1]$ set iteratively $x_i:= w_{\Sigma_{u[i]}} \cdot v_i \cdot u[i]$ for $v_i=\epsilon$ if $u[i]=u[i+1]$ and $v=u[i+1]$ otherwise. Lastly, set $w= x_1 \cdots x_{k-1} \cdot w_{\Sigma_{u[k]}}$.

Firstly, we want to prove that $u$ is a absent scattered factor of the returned word $w$. 
By the construction of $x_i$, we get $\ar_i(w)=x_i$: the prefix of $x_i$ of length $\sigma-1$ contains all letters
of $\Sigma$ but $u[i]$, then $u[i+1]$ is appended iff $u[i]\neq u[i+1]$ and lastly $u[i]$ is appended which is therefore unique in $x_i$. Since $w_{\Sigma_{u[k]}}$ contains all letters of $\Sigma$ but $u[k]$, we get $u[1..k-1]=\m(w)$ and  $u\not\in\ScatFact_k(w)$. To prove $w\in\NUniv_{\Sigma,1,k}$ we show the three conditions of Theorem~\ref{niceequivalence}. We already showed $\iota(w)=k-1$. Moreover, $\alph(w_{\Sigma_{u[k]}})=\Sigma\backslash\{u[k]\}$ follows by definition. Set $y_i=v_iu[i]w_{\Sigma_{u[i+1]}}$ for all $i\in[k-1]$. By the definition of $v_i$
we get that $y_{k-1}^Ry_{k-2}^R\cdots y_1^R w_{\Sigma_{u[1]}}^R$ is the arch factorisation of $w^R$. This implies
$\alph(w^R)=\Sigma\backslash\{u[1]\}$ and thus the second condition is fulfilled. Hence, only the third conditions remains to be proven. \\
\textbf{case 1:} $k$ even\\
Set 
\begin{align*}
u_1&=x_1\cdots x_{\frac{k}{2}},& v_1&=y_{k-1}\cdots y_{\frac{k}{2}+1}, \\
u_2&=x_1\dots x_{\frac{k}{2}-1},& v_2&=y_{k-1}\cdots y_{\frac{k}{2}}.
\end{align*}
Thus, we get $x_1=w_{\Sigma_{u[\frac{k}{2}+1]}}$ and $x_2=w_{\Sigma_{u[\frac{k}{2}]}}$. Both fulfil by definition the required property.\\
\textbf{case 2:} $k$ odd\\
Set 
\[
u=x_1\cdots x_{\frac{k-1}{2}}\mbox{ and } v=y_{k-1}\cdots y_{\frac{k+1}{2}}.
\]
Thus, we get $x=w_{\Sigma_{u[\frac{k-1}{2}+1]}}$ which fulfils by definition the required property.\\
Hence, in both cases all three properties are fulfilled and by Theorem~\ref{niceequivalence}, we have $w\in\NUniv_{\Sigma,1,k}$.

Assuming that all $w_{\Sigma_\ta}$ for all $\ta\in\Sigma$ are precalculated, we just have to compare $u[i]$ with $u[i+1]$
and append the appropriate words for obtaining $w$. Since $u\in\Sigma^k$, we have $k$ of those comparisons and extensions
of the word.

It remains to show that $w$ is of minimal length among all nearly $k$-universal words where $u$ is the only absent scattered factor. Let $w'\in\NUniv_{\Sigma,1,k}\backslash\{w\}$ with $u\not\in\ScatFact_k(w')$. By Theorem~\ref{theorem_univ} we know $\iota(w')=k-1$ and by Corollary~\ref{absent} we have $\m(w')\ta_{w'}=u$. This implies
immediately $\m(w)=\m(w')$ and $\ta_{w'}=\ta_w$. Since each arch has to contain the complete alphabet, $w'$ has at least 
one arch which is shorter than a corresponding arch in $w$. By the definition of $w$ we know that this arch $\alpha$ contains each letter of $\Sigma$ exactly once. Thus, in  the arch factorisation of $(w')^R$ the changed arch goes further to the left. The corresponding $x$ from Theorem~\ref{niceequivalence} does not contain all letters from $\Sigma$ but one and we can conclude that $w'$ is not nearly $k$-universal.

This concludes the proof.\qed	
\end{proof}

\fi

\ifpaper
\else
\SetKwComment{Comment}{/* }{ */}
\begin{algorithm}\label{alg:one}
	\caption{Computing $w\in\NUniv_{\Sigma,1,k}$ for $u\in\Sigma^k$ absent (cf. Theorem~\ref{utow}).}
	\KwData{Given $u \in \Sigma^k$ with $\Sigma = \{\ta_1, \dots \ta_\sigma\}$.}
	\KwResult{nearly $k$-universal word $w \in \Sigma^*$ with $\operatorname{ScatFact}_k(w) = \Sigma^k \setminus \{u\}$}
	$w := \epsilon$\;
	$w_\Sigma = \ta_1 \cdots \ta_\sigma$\;
	\For{$i = 1$ to $k-1$}{
		\eIf{$u[i] \neq u[i+1]$}{
			$w := w \cdot w_{\Sigma_{u[i]}} \cdot u[i+1] \cdot u[i]$
		}{$w := w \cdot w_{\Sigma_{u[i]}} \cdot u[i]$}
	}
	$w := w \cdot w_{\Sigma_{u[k]}}$\;
	\Return{$w$}\;
\end{algorithm}

\begin{figure}
\begin{tikzpicture}[scale=0.5,
blackbox/.style={rectangle, draw=black, thick, minimum height=0.6cm},
redbox/.style={rectangle, draw=black!75, 
			 thick, minimum height=0.5cm}
]

\draw[-] (0,0.75) -- (0,1.25);
\draw[-] (24,0.75) -- (24,1.25);


\coordinate (1) at (0,1);
\coordinate (2) at (4,1);
\coordinate (3) at (8,1);
\coordinate (4) at (12,1);
\coordinate (5) at (16,1);
\coordinate (6) at (20,1);
\coordinate (7) at (24,1);
\coordinate (r1) at (24,1.5);
\coordinate (r2) at (24,0.5);

\draw[-|]  (1) -- (2) {};
\draw[-|]  (2) -- (3) {};
\draw[-|]  (3) -- (4) {};
\draw[-|]  (4) -- (5) {};
\draw[-|]  (5) -- (6) {};
\draw[-|]  (6) -- (7) {};

\node (m1) at ($(2) + (-0.5, 0.5)$) {$\color{gray} \ta$};  
\node (m2) at ($(3) + (-0.5, 0.5)$) {$\color{gray} \tb$};  
\node (m3) at ($(4) + (-0.5, 0.5)$) {$\color{gray} \tc$}; 
\node (m4) at ($(5) + (-0.5, 0.5)$) {$\color{gray} \tc$};  
\node (m5) at ($(6) + (-0.5, 0.5)$) {$\color{gray} \ta$};  


\draw [->,gray] (1.south west) to [out=-30,in=-150] (2.south);
\draw [->,gray] (2.south) to [out=-30,in=-150] (3.south);
\draw [->,gray] (3.south) to [out=-30,in=-150] (4.south);
\draw [->,gray] (4.south) to [out=-30,in=-150] (5.south);
\draw [->,gray] (5.south) to [out=-30,in=-150] (6.south);
\draw [->,gray] (6.south) to [out=-30,in=-170] (r2);

\node[redbox] (i1) at ($(2) + (-1.5,0.5)$) {$\tb$};
\node[redbox] (i2) at ($(3) + (-1.5,0.5)$) {$\tc$};
\node[redbox] (i3) at (m3) {$\color{gray} \tc$};
\node[redbox] (i4) at ($(5) + (-1.5,0.5)$) {$\ta$};
\node[redbox] (i5) at ($(6) + (-1.5,0.5)$) {$\tb$};

\draw [->,gray, dotted] (r1.north) to [out=150,in=30] (i5.north west);
\draw [->,gray, dotted] (i5.north west) to [out=150,in=30] (i4.north west);
\draw [->,gray, dotted] (i4.north west) to [out=150,in=30] (i3.north west);
\draw [->,gray, dotted] (i3.north west) to [out=150,in=30] (i2.north west);
\draw [->,gray, dotted] (i2.north west) to [out=150,in=30] (i1.north west);
\draw [->,gray, dotted] (i1.north west) to [out=150,in=30] ($(1) + (0,1)$);

\end{tikzpicture}
\caption{An illustration for the construction for the absent scattered factor $u = 
\mathtt{abccab}$.}\label{newlabel}
\end{figure}
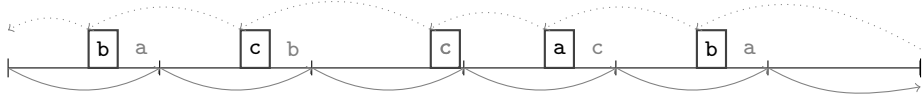

\fi


As illustrated in Figure~\ref{newlabel}, let $u=\ta\tb\tc\tc\ta\tb$ and $\bullet$ represent placeholder. Since $u[1..5]$ is $\m(w)$, we get 
$(\bullet\ta)\cdot(\bullet\tb)\cdot(\bullet\tc)\cdot(\bullet\tc)\cdot(\bullet\ta)$. By 
$\alph(\r(w))=\Sigma\backslash\{\tb\}$, we get 
$(\bullet\ta)\cdot(\bullet\tb)\cdot(\bullet\tc)\cdot(\bullet\tc)\cdot(\bullet\ta)\cdot\ta\tc$. 
Including the arches of $w^R$ we obtain 
$(\bullet\tb\ta)\cdot(\bullet\tc\tb)\cdot(\bullet\tc)\cdot(\bullet\ta\tc)\cdot(\bullet\tb\ta)
\cdot\ta\tc$. Now, the $\bullet$ are replaced by the missing letters from each arch of $w$. Thus, 
we finaly get  
$(\tb\tc\tb\ta)\cdot(\ta\tc\tc\tb)\cdot(\ta\tb\tc)\cdot(\ta\tb\ta\tc)\cdot(\tb\tc\tb\ta)
\cdot\ta\tc$.

\begin{remark}\label{length}
Notice that the length of the resulting nearly $k$-universal word $w$ depends on the given absent 
scattered factor $u \in \Sigma^k$. If $\cond(u)=u_1\cdots u_r$ for an $r\in\N$, we have 
$|w|=k\sigma+r-2$. Thus, if $u$ is unary, we have 
$|w|=k\sigma-1$.
\end{remark}

\begin{definition}
Let $w_u\in\NUniv_{\Sigma,1,k}$ be the unique word of minimal length w.r.t. a given order 
$<_{\Sigma}$ with $u\not\in\ScatFact_k(w_u)$.
\end{definition}

\begin{corollary}\label{nearlysimon}
Given $k\in\N$, we have $|\NUniv_{\Sigma,1,k}/\sim_k|=\sigma^k$, i.e., restricting the Simon 
congruence to nearly $k$-universal words leads to $\sigma^k$ different congruence classes.
\end{corollary}

\ifpaper
\else
\begin{proof}
Follows directly with Theorem~\ref{utow}.\qed
\end{proof}

\fi

By Corollary~\ref{nearlysimon} we know how many congruence classes in $\NUniv_{\Sigma,1,k}$ 
w.r.t. $\sim_k$ exist. Now we show when $w\sim_kw_u$ holds for $w\in\Sigma^{\ast}$,
i.e., we characterise $[w_u]_{\sim_k}$. Therefore, we need some further insights into nearly 
$k$-universal words.

\begin{lemma}{} \label{paperstream}
	Given $w \in \NUniv_{\Sigma,1,k}$, we have $\ar_1(w) 
	\cdots \ar_{i-1}(w)\alpha\operatorname{ar}_{i+1}(w) \cdots$ $\operatorname{ar}_{k-1}(w) \beta \in\NUniv_{\Sigma,1,k}$
	for all $i\in[k-1]$, if
	$\alpha[|\alpha|]=\m(w)[i]$, $\alph(\alpha[1..|\alpha|-1])=\alph(\inner_i(w))$, $\inner_i(w)\in\ScatFact(\alpha[1..|\alpha|-1])$, $|\r(w)|\leq|\beta|$, and $\alph(\beta)=\alph(\r(w))$.
\end{lemma}
\ifpaper
\else    
\begin{proof}
	Let $w\in\NUniv_{\Sigma,1,k}$. By Theorem \ref{theorem_univ} we know that exactly $\m(w)\ta_w$ is the only absent scattered factor. By the conditions, neither the modus nor the alphabet of the rest are changed, i.e. $\m(w)\ta_w$
	is still absent. Since adding letters which are not the modus, does not change the modus, neither deletes scattered factors we know that $\m(w)\ta_w$ is still the only absent scattered factor.\qed
\end{proof}

\fi

Let $P(w)$ be the set 
of all words obtainable from $w$ by Lemma~\ref{paperstream}.

\begin{remark}
Lemma~\ref{paperstream} implies that for all $n\in\N$ and $u\in\Sigma^k$ with $n\geq |w_u|$ there 
exists $w\in\NUniv_{\Sigma,1,k}\cap\Sigma^n$, i.e. $|[w_u]_{\sim_k}|=\infty$. 
\end{remark}

%
%

We are now able to give a characterisation of the congruence classes of $\sim_k$ in 
$\NUniv_{\Sigma,1,k}$. Since we know that for each $u\in\Sigma^k$ there
exists one congruence class, we fix $u\in\Sigma^k$. We know so far that for each $w$ obtained by the application of Lemma~\ref{paperstream}, we have $w\in[w_u]_{\sim_k}$. Notice that Lemma~\ref{paperstream} cannot be generalised to an equivalence, since deleting letters from arches may violate
the nearly $k$-universality: considering $w = (\ta\ta\tb)\cdot \tb$ and deleting one $\ta$ 
in the first arch, indeed does not change
the modus, but it deletes $\ta\ta$ and therefore we have $|\ScatFact_2(\ta\tb\tb)|<3$. Recall that 
the output $w_u$ of the algorithm in Theorem~\ref{utow} is w.r.t. a given order $<_{\Sigma}$ on 
$\Sigma$, in particular $\Pref_{\sigma-1}(\ar_i(w_u))$, for all $i\in[k-1]$, is the 
lexicographically smallest word containing all letters of $\Sigma$ but $\m[i]$. Analogously, 
$\r(w_u)$ is the lexicographically smallest word
containing all letters but $u[k]$. If we change this order, we obtain other words of the same length, which are all by Theorem~\ref{utow} of minimal length. Moreover, if we choose different orders for each arch and for the rest, we still obtain a nearly $k$-universal word since
the crucial point of Theorem~\ref{utow} still holds. Thus, each such word can be obtained from $w_u$ 
by applying some morphic permutation of
$\Sigma$ on $\Pref_{\sigma-1}(\ar_i(w_u))$ and $\r(w_u)$ for all $i\in[k-1]$.

\begin{definition}
Let $\pi_1,\ldots,\pi_{\sigma!}$ be the different morphic permutations on $\Sigma$,
set $p_i=\Pref_{\sigma-1}(\ar_i(w_u))$ for all $i\in[k-1]$, and choose $s_1,\dots,s_{k-1}\in\Sigma^{\ast}$ with $w_u=p_1s_1\cdots p_{k-1}s_{k-1}\r(w_u)$.
Define the {\em basis} of $[w_u]_{\sim_k}$ by $
B_u=\{w\in\Sigma^{\ast}|\,\exists i_1,\dots,i_k\in[\sigma!]: 
w=\pi_{i_1}(p_1)s_1\cdots \pi_{i_{k-1}}(p_{k-1})s_{k-1}\pi_{i_k}(\r(w_u))\}$.
\end{definition}

\begin{remark}
For $u\in\Sigma^k$, we have $|B_u|=((\sigma-1)!)^{k-1}(\sigma-1)!$.
\end{remark}

Based on this $B_u$ and Lemma~\ref{paperstream}, we can characterise $[w_u]_{\sim_k}$.

\begin{theorem}\label{charclass}
Given $u\in\Sigma^k$, we have $[w_u]_{\sim_k}=\{w\in\Sigma^{\ast}|\,\exists v\in B_u:w\in P(v) \}$.
\end{theorem}

\ifpaper
\else
\begin{proof}
If $w\in P(v)$ for some $v\in B_u$, we have immediatly $w\in[w_u]_{\sim_k}$. Assume $w\in[w_u]_{\sim_k}$. Thus $w\in\NUniv_{\Sigma,1,k}$
and $\m(w)\ta_w=u$. Now, we examine $w$'s \nth{i} arch for a fixed $i\in[k-1]$. We know $\alph(\inner_i(w))=\Sigma\backslash\{u[i]\}$. Let
$u[i]\neq u[i+1]$. Suppse that $|\inner_i(w)|_{u[i+1]}=1$. The application of Theorem~\ref{niceequivalence} with $\hat{k}=i-1$
and $\tilde{k}=k-i$ implies that the \nth{$(k-i)$}-arch from $v^R$ ends in this occurrence of $u[i+1]$, i.e. $u[i],u[i+1]\not\in\alph(x)$.
Since this is a contradiction to $w\in\NUniv_{\Sigma,1,k}$ we not only have $|\inner_i(w)|_{u[i+1]}\geq 2$ but also Theorem~\ref{niceequivalence} leads to $\ar_i(w)=\alpha_iu[i+1]\beta_iu[i]$ with $\alph(\alpha_i)=\Sigma\backslash\{u[i]\}$ and
$\beta_i\in(\Sigma\backslash\{u[i]\})^{\ast}$. Thus, there exists $v_i\in\ScatFact_{\sigma-1}(\alpha_i)$ with $\alph(v)=\Sigma\backslash\{u[i]\}$.
Hence there exists a permutation $\pi$ on $\Sigma$ which morphically applied yields $\pi(v_i)=\Pref_{\sigma-1}(\ar_i(w_u))$.
Since $\alph(\r(w))=\Sigma\backslash\{u[k]\}$ we get by the same argument an $r$ which is a permuation of $\r(w_u)$. This leads
to $v=v_1u[1]\cdots v[k-1]u[k-1]r\in B_u$. Adding all letters of $\alpha_i,\beta_i$ and $\r(w)$, resp., which are not in $v_i$ and $r$, resp.,
implies $w\in P(v)$.\qed
\end{proof}

\fi

Let $u=\ta\tb\tb\tc$ and $\ta<\tb<\tc$. By Theorem~\ref{utow} we get $w_u=(\tb\tc\tb\ta)\cdot(\ta\tc\tb)\cdot(\ta\tc\tc\tb)\cdot\ta\tb$
and $w\in B_u$ iff $w=w_1w_2w_3w_4$ with $w_1\in\{\tb\tc\tb\ta,\tc\tb\tb\ta\}$, $w_2\in\{\ta\tc\tb,\tc\ta\tb\}$, $w_3\in\{\ta\tc\tc\tb,
\tc\ta\tc\tb\}$, $w_4\in\{\ta\tb,\tb\ta\}$. Thus, we have $16$ basis elements for $u$. Each this 
word can be enriched by additional letters in the inner of an arch and the rest w.r.t. 
Lemma~\ref{paperstream} to obtain all elements equivalent to $w_u$.

We finish this section with a third characterisation of nearly $k$-universal words that relies on  
 Theorem~\ref{niceequivalence} and Lemma~\ref{cutting_first} and illustrates the relation of $w$ 
and $w^R$ in $\NUniv_{\Sigma,1,k}$.

\begin{theorem} {} \label{recursiv}
We have $w\in\NUniv_{\Sigma,1,k}$ iff $\iota(w)=k-1$, $|\alph(\r(w))|=\sigma-1$, and $(\operatorname{ar}_2(w^R)$ ... $\operatorname{ar}_{k-1}(w^R)\r(w^R))^R\in\NUniv_{\Sigma,1,k-1}$.
\end{theorem}
\ifpaper
\else
\begin{proof}
	Consider first $w\in\NUniv_{\Sigma,1,k}$. The first two conditions follow by Theorem~\ref{theorem_univ}. Since $w\in\NUniv_{\Sigma,1,k}$,
	we have $w^R\in\NUniv_{\Sigma,1,k}$ and Lemma~\ref{cutting_first} implies $\hat{w} = \operatorname{ar}_2(w^R)
	\cdots\operatorname{ar}_{k-1}(w^R)\operatorname{r}(w^R)\in\NUniv_{\Sigma,1,k}$. Thus, we have $\hat{w}^R\in\NUniv_{\Sigma,1,k-1}$.
	
	Consider now $w\in\Sigma^{\ast}$ with $\iota(w)=k-1$, $|\alph(\r(w))|=\sigma-1$, and
	$\hat{w} = (\operatorname{ar}_2(w^R)
	\cdots\operatorname{ar}_{k-1}(w^R)\operatorname{r}(w^R))^R\in\NUniv_{\Sigma,1,k}$. Thus, we have
	$\ar_i(\hat{w})=\ar_i(w)$ for all $i\in[k-2]$. By $\hat{w}\in\NUniv_{\Sigma,1,k-2}$ we get
	$|\alph(\r(\hat{w}))|=\sigma-1$. Thus, $w=\hat{w}\ar_1(w^R)$ fulfils the conditions of Theorem~\ref{niceequivalence} and the claim is proven.\qed
\end{proof}

\fi

Notice that only the deletion of a reversed arch from the beginning leads to an equivalence. 
Deleting the first arch of $w$ does not suffice for a characterisation as witnessed by 
$w=\tb\tc\ta\ta\tb\tc\ta\tb$: indeed, we have $\iota(\ta\tb\tc\ta\tb)=1$, 
$\alph(\r(w))=\Sigma\backslash\{\tc\}$, and $\mathtt{abcabcab} \in \NUniv_{\Sigma,1,2}$
but we get  $\mathtt{(bca) \cdot (abc) \cdot ab}\not \in \NUniv_{\Sigma,1,3}$.

In this section, we presented a characterisation for nearly $k$-universal words as well as
the index of $\sim_k$ and a characterisation of its congruence classes.

\section{$m$-Nearly $k$-Universal Words}\label{mnearly}
\label{sec:mnearly}
In this section, we consider $m$-nearly $k$-universal words, where $m$ is not necessarily $1$, i.e., we
are interested in $w\in\Sigma^{\ast}$ with $|\ScatFact_k(w)|=\sigma^k-m$.
Implicitly, a subset of these words was investigated in \cite{kosche2021absent}. There, the authors 
determine all shortest absent scattered factors, i.e. if $\iota(w)=k-1$ and 
$|\ScatFact_k(w)|=\sigma^k-m$, we have that $w\in\NUniv_{\Sigma,m,k}$.
In contrast to $1$-nearly $k$-universal words, for $m>1$, $\iota(w)=k-1$ does not necessarily hold 
as witnessed by $\mathtt{ababca} \in \NUniv_{\Sigma,14,3}$ with $\iota(\mathtt{ababca}) = 1 \neq 
2$. Thus, a thorough characterisation of $\NUniv_{\Sigma,m,k}$ is 
still open. Unfortunately, we cannot give such a characterisation but we present some first 
insights for $m\in\{\sigma^k,\sigma^{k-1},2\}$ as well as a full characterisation of the subclass 
established in \cite{kosche2021absent} including the congruence classes of $\sim_k$ in this case.

\begin{remark}\label{nothing}
Similar to $\NUniv_{\Sigma,0,k}=\Univ_{\Sigma,k}$, the set $\NUniv_{\Sigma,\sigma^k,k}$ provides
exactly one equivalence class for $\sim_{k}$, since exactly the words strictly shorter than $k$ do not have
any scattered factor of length $k$. 
\end{remark}

Now, we have a look at $m\in\{\sigma^k-1,\sigma^k-2\}$. Since 
$w\in \NUniv_{\Sigma,\sigma^k-1,k}$ for all $w\in\Sigma^k$, we have $|u|\geq k+1$ 
for all $u\in\NUniv_{\Sigma,m,k}$ with $m<\sigma^k-1$.

\begin{proposition}\label{sigmak-1}
For each $k\in\N$, we have $|\NUniv_{\Sigma,\sigma^k-1,k}/\sim_k|=\sigma^k$.
\end{proposition}

\ifpaper
\else
\begin{proof}
	 First, we can observe that for $u,v\in\Sigma^k$ with $u\neq v$, we have $u,v\in\NUniv_{\Sigma,\sigma^k-1,k}$
	 and $[u]_{\sim_k}\neq[v]_{\sim_k}$. Now let $w\in\Sigma^{\ast}$ with $|u|\geq k$. If $\alph(w)=\{\ta\}\subset\Sigma$,
	 we have immediately that all scattered factors but $\ta^k$ are absent. Thus, we have $[\ta^k]_{\sim_k}=\{\ta\}^{\ast}$
	 for all $\ta\in\Sigma$. If, on the other hand, we have $|\alph(w)|\geq 2$, we can factorise $w=w_1\ta w_2\tb w_3$ with
	 $\ta\not\in\alph(w_1)$ and $\tb\not\in\alph(w_3)$ and obtain that $w[1..k]$ and $w[|w|-k+1..|w|]$ are different scattered factors of $w$.
	 Thus, all $w\in\NUniv_{\Sigma,\sigma^k-1,k}$ are unary. This proves the claim.\qed
\end{proof}

\fi

\begin{lemma}\label{2insidek+1}
If $w\in\NUniv_{\Sigma,\sigma^k-2,k}$ then $|\alph(w)|=2=|\cond(w)|$.
\end{lemma}

\ifpaper
\else
\begin{proof}
If $|\alph(w)|\leq 1$, we have $|\ScatFact_k(w)|\in\{0,1\}$.
Suppose $|\alph(w)|\geq 3$. Then there exists $\ta_1,\ta_2,\ta_3\in\Sigma$. Choose words  $w_1,w_2,w_3,w_4\in\Sigma^{\ast}$ and $r_1,r_2,r_3\in\N$
such that
\[
w=w_1\ta_1^{r_1}w_2\ta_2^{r_2}w_3\ta_3^{r_3}w_4
\]
and $w_1[|w_1|],w_2[1]\neq\ta_1$, $w_2[|w_2|],w_3[1]\neq\ta_2$, and $w_3[|w_3|],w_4[1]\neq\ta_3$. Moreover, choose
$v_1\in\Pref(w_1\ta_1^{r_1-1})$, $v_2\in\Pref(w_2\ta_2^{r_2-1})$, $v_3\in\Pref(w_3\ta_3^{r_3-1})$, and $v_4\in\Pref(w_4)$
such that $|v_1|+|v_2|+|v_3|+|v_4|=k-2$. Set 
\[
u_1= v_1  v_2 \ta_2 v_3\ta_3 v_4,\, u_2= v_1 \ta_1 v_2   v_3\ta_3v_4,\, u_3= v_1  \ta_1 v_2 \ta_2 v_3v_4.
\]
Then we have $u_1,u_2,u_3\in \ScatFact_k(w)$. With $\ell_1=|v_1|+1$.
By $u_1[\ell_1]=v_2[1]\neq\ta_1= u_2[\ell_1],u_3[\ell_1]$, we get $u_1\neq u_2,u_3$. 
Let $\ta_3^{s}\in\Suff(v_3)$ maximal for $s\in\N_0$. Then $\ta_3^{s+1}v_4\in\Suff(u_2)$ but $\ta_3^{s+1}v_4\not\in\Suff(u_3)$.
Thus, $u_2\neq u_3$ - a contradiction to $w\in\NUniv_{\Sigma,\sigma^k-2,k}$.

\medskip

Suppose 
there exist $r_1,r_2,r_3\in\N$ and
$w_1\in\Sigma^{\ast}$ with $w=\ta^{r_1}\tb^{r_2}\ta^{r_3}w_1$, i.e., $|\cond(w)|>2$ (we assume 
w.l.o.g. that $w[1]=\ta$).
Set $u_1=\ta^{s_1-1}\tb^{s_2}\ta^{s_3}v_1$, $u_2=\ta^{s_1}\tb^{s_2-1}\ta^{s_3}v_1$, 
$u_3=\ta^{s_1}\tb^{s_2}\ta^{s_3-1}v_1$ such that $s_1-1+s_2+s_3+|v_1|=k$
for $s_i\leq r_i$, $i\in[3]$, $v_1\in\Pref(w_1)$. Then $u_1,u_2$, and $u_3$ are different scattered 
factors of $w$ - a contradiction.\qed

\end{proof}

\fi

\begin{proposition}\label{sigmak-2}
	For each $k\in\N$, we have $|\NUniv_{\Sigma,\sigma^k-2,k}/\sim_k|= 2 \binom{\sigma}{2}(k+2)$.
\end{proposition}

\ifpaper
\else
\begin{proof}
Let $w\in\NUniv_{\Sigma,\sigma^k-2,k}$. By Lemma~\ref{2insidek+1} we have $|\alph(w)|=2$ and there 
exist $r_1,r_2\in\N$
such that w.l.o.g. $w=\ta^{r_1}\tb^{r_2}$ with $r_1+r_2\geq k+1$. Thus all $\ta^{s_1}\tb^{s_2}$ are scattered factors of $w$ of length $k$
with $s_1+s_2=k$, $s_1\leq r_1, s_2\leq r_2$. \\
\textbf{case 1:} $r_1+r_2=k+1$\\
In this case, we have exactly two scattered factors, namely $\ta^{r_1-1}\tb^{r_2}$ and
$\ta^{r_1}\tb^{r_2-1}$ (choosing less $\ta$ required more $\tb$ than available and v.v.).
Thus, for fixed $\ta_1,\ta_2\in\Sigma$, all classes $[\ta_1^{t_1}\ta_2^{t_2}]_{\sim_k}$
are different, for all possible $t_1,t_2\in\N$ with $t_1+t_2=k+1$. These are $2k$ congruence classes for each choice of $\ta_1,\ta_2$. Thus, we have $2k\binom{\sigma}{2}$ classes. Notice that for all these congruence classes $[\ta_1^{t_1}\ta_2^{t_2}]$, we have $|\alph(u)|=2$ for all $u\in\ScatFact_k(\ta_1^{t_1}\ta_2^{t_2})$.\\
\textbf{case 2:} $r_1+r_2>k+1$, $r_1,r_2<k$\\
In this case, we have $\min\{r_1,r_2\}$ different scattered factors and by the choice
of $w$, we have $\min\{r_1,r_2\}=2$ leading to $w\in\{\ta^2\tb^{r_2},\ta^{r_1}\tb^2\}$
with $\ta^2\tb^{k-2}, \ta\tb^{k-1}$ and $\ta^{k-2}\tb$, $\ta^{k-1}\tb$, resp., as scattered factors. By $\ta^2\tb^{r_2}\in[\ta^2\tb^{k-1}]_{\sim_k}$ and $\ta^{r_1}\tb^2\in[\ta^{k-1}\tb^2]_{\sim_k}$, we do not obtain different classes.\\
\textbf{case 3:} $r_1+r_2>k+1$ and $r_1\geq k$ or $r_2\geq k$\\
In this case, we have $|\ScatFact_k(w)|=\min\{r_1,r_2\}+1$. By the choice of $w$,
we have $2=\min\{r_1,r_2\}+1$, thus $1=\min\{r_1,r_2\}$ leading to $w\in\{\ta\tb^{r_2},\ta^{r_1}\tb\}$ with the scattered factors $\ta\tb^{k-1}$, $\tb^{k}$ and $\ta^{k-1}\tb$, $\ta^k$ resp. This implies that in this case we get new congruence classes, since all
words have one unary and one binary scattered factor. Thus, for fixed $\ta_1,\ta_2\in\Sigma$, the four classes $[\ta_1\ta_2^t]_{\sim_k}$, $[\ta_1^t\ta_2]_{\sim_k}$,
$[\ta_2\ta_1^t]_{\sim_k}$, and $[\ta_2^t\ta_1]_{\sim_k}$, for a fixed $t>k$ are different. Thus, we have $4\binom{\sigma}{2}$ new classes.\\
Summing up, we get $|\NUniv_{\Sigma,\sigma^k-2,k}/\sim_k|=2\binom{\sigma}{2}(k+2)$.\qed
\end{proof}

\fi

Proposition~\ref{sigmak-2} shows that the formula determining the index of $\sim_k$ gets more 
complicated the farther $m$ is
from $0$ or $\sigma^k$, resp. Now we show a similar result to Theorem~\ref{theorem_univ} for 
$\NUniv_{\Sigma,2,k}$ backing the observation that the conditions on $w$ get more complicated. 
Notice that Theorem~\ref{m2} does not hold for $\sigma=2$ witnessed by 
$w=\ta\ta\tb\ta\ta\in\NUniv_{\Sigma,4,3}$ but $\iota(w)=1$. Moreover,
$w\in\NUniv_{\Sigma,\sigma-1,k}$ implies $\iota(w)=k-1$.

\begin{theorem}\label{m2}
Let $w\in\NUniv_{\Sigma,2,k}$ with $\sigma>2$. Then $\iota(w)=k-1$ and either 
$|\alph(\r(w))|=|\alph(\r(w))|=\sigma-1$, or $|\alph(\r(u))|=\sigma-1$ and 
$|\alph(\r(u^R))|=\sigma-2$ for all $u\in\{w,w^R\}$.
\end{theorem}

\ifpaper
\else
\begin{proof}
Suppose $\iota(w)<k-1$. Then there exits $v\in\Sigma^{k-1}$ with $v\not\in\ScatFact_k(w)$.
Since $\sigma>2$ there exist $\ta_1,\ta_2,\ta_3\in\Sigma$. This implies 
$v\ta_i\not\in\ScatFact_k(w)$ for all $i\in[3]$. Thus, we know $\iota(w)=k-1$.

If we had $w \in \PerfUniv_{\Sigma,k-1}$, each $m(w)\ta$ for $\ta\in\Sigma$ would be an absent 
scattered factor. Thus, we have $r(w)\neq\varepsilon$. Analogously we get $r(w^R)\neq\varepsilon$. 
Similarly to the previous argumentation,
if $|\alph(\r(w))|<\sigma-2$ or $|\alph(\r(w^R))|<\sigma-2$, we would have at least three absent 
scattered factors. In the case of $\alph(\r(w)) = \Sigma \setminus \{\ta,\tb\}$ for $\ta,\tb \in 
\Sigma$ with $\ta \neq \tb$ we know that $\m(w)\ta, \m(w)\tb \notin \ScatFact_k(w)$.
Thus, $\ta \m(w)^R, \tb \m(w)^R \notin \ScatFact_k(w^R)$. 
Thus, $\ta,\tb \notin \alph(r(w^R))$. If we had $\m(w)^R\neq \m(w^R)$, then $\m(w^R)^R\ta$ would be 
a third absent scattered factor. Thus, we have $m(w)^R=m(w^R)$. This leads to the following 
contradiction: by $\ar_1(w)=\r(w^R)^R\m(w^R)[k-1]=\r(w^R)^R\m(w)[1]$ and 
$\ta,\tb\not\in\alph(\r(w^R))$ we would get that either $\ta$ or $\tb$ cannot be in $\ar_1(w)$.\qed
\end{proof}

\fi

We finish this section by characterising $\NUniv_{\Sigma,m,k}\cap\Univ_{\Sigma,k-1}$.
Let from now on $w\in\Univ_{\Sigma,k-1}$.
By $|\alph(\r(w^R))| <\sigma$, we have $\r(w)^R\in\Pref(\ar_1(w^R))$ and 
$\m(w^R)[1]\in\alph(\ar_{k-1}(w))$. Thus, choose $\alpha_{k-1},\beta_{k-1}\in\Sigma^{\ast}$ 
with $\ar_{k-1}(w) = \alpha_{k-1}\beta_{k-1}$ and $\ar_1(w^R) = (\beta_{k-1} \r(w))^R$.
With $\alph(\beta_{k-i}) \subseteq \Sigma$, inductively there exist $\alpha_i, 
\beta_i \in\Sigma^\ast$ such that $\ar_{k-i}(w) = \alpha_i \beta_i$ and $\ar_i(w^R) = 
(\beta_i \alpha_{i+1})^R$ with $\alpha_{k} = \r(w)$ and $\alpha_1 = \r(w^R)^R$, for all $i\in[k-1]$.

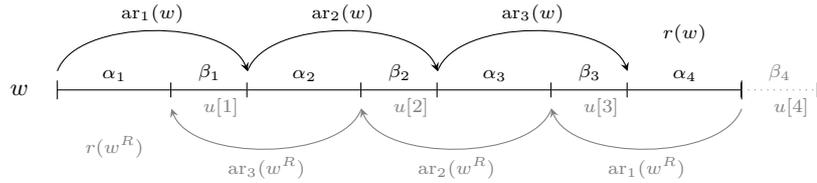
\begin{figure}
 \centering
\begin{tikzpicture}[scale=0.5,
blackbox/.style={rectangle, draw=black, thick, minimum height=0.6cm},
redbox/.style={rectangle, draw=red!75, 
			 thick, minimum height=0.5cm}
]

\coordinate (start) at (0,1);
\coordinate (end) at (18,1);

\coordinate (ar1) at (5,1);
\coordinate (ar2) at (10,1);
\coordinate (ar3) at (15,1);
\coordinate (r) at (18,1);

\coordinate (rR) at (0,1);
\coordinate (ar1R) at (13,1);
\coordinate (ar2R) at (8,1);
\coordinate (ar3R) at (3,1);

\node[black] (w) at (-1,1) {$w$};
\draw[-,black] (0,0.75) -- (0,1.25);
\draw[-,black] (18,0.75) -- (18,1.25);

\draw[black, -stealth] ($ (start)+(0,0.5) $) to[out=70, in=110, distance=1.5cm] node[above, black, font=\scriptsize] {$\ar_1(w)$} ($ (ar1)+(0,0.5) $);
\draw[black, -stealth] ($ (ar1)+(0,0.5) $) to[out=70, in=110, distance=1.5cm] node[above, black, font=\scriptsize] {$\ar_2(w)$} ($ (ar2)+(0,0.5) $);
\draw[black, -stealth] ($ (ar2)+(0,0.5) $) to[out=70, in=110, distance=1.5cm] node[above, black, font=\scriptsize] {$\ar_3(w)$} ($ (ar3)+(0,0.5) $);
\node[black, font=\scriptsize] (r) at(16.5, 2.5) {$r(w)$};

\draw[gray, -stealth] ($ (end)+(0,-0.5) $) to[out=-110,in=-70, distance=1.5cm] node[below, gray, font=\scriptsize] {$\ar_1(w^R)$}  ($ (ar1R)+(0,-0.5) $);
\draw[gray, -stealth] ($ (ar1R)+(0,-0.5) $) to[out=-110,in=-70, distance=1.5cm] node[below, gray, font=\scriptsize] {$\ar_2(w^R)$} ($ (ar2R)+(0,-0.5) $);
\draw[gray, -stealth] ($ (ar2R)+(0,-0.5) $) to[out=-110,in=-70, distance=1.5cm] node[below, gray, font=\scriptsize] {$\ar_3(w^R)$} ($ (ar3R)+(0,-0.5) $);
\node[gray, ->, font=\scriptsize] (r) at(1.5,-0.5) {$r(w^R)$};

\draw[-|,black] (0,1) -- node[above,font=\scriptsize] {$\alpha_1$} (ar3R);
\draw[-|,black] (ar3R) -- node[above,font=\scriptsize] {$\beta_1$} (ar1);
\draw[-|,black] (ar1) -- node[above,font=\scriptsize] {$\alpha_2$} (ar2R);
\draw[-|,black] (ar2R) -- node[above,font=\scriptsize] {$\beta_2$} (ar2);
\draw[-|,black] (ar2) -- node[above,font=\scriptsize] {$\alpha_3$} (ar1R);
\draw[-|,black] (ar1R) -- node[above,font=\scriptsize] {$\beta_3$} (ar3);
\draw[-|,black] (ar3) -- node[above,font=\scriptsize] {$\alpha_4$} (end);
\draw[-|,gray, dotted] (18,1) -- node[above,font=\scriptsize] {$\beta_4$} (20,1);

\node[gray,font=\scriptsize, below left] at (ar1) {$u[1]$};
\node[gray,font=\scriptsize, below left] at (ar2) {$u[2]$};
\node[gray,font=\scriptsize, below left] at (ar3) {$u[3]$};
\node[gray,font=\scriptsize, below left] at (20,1) {$u[4]$};

\end{tikzpicture}
\caption{$\alpha$-$\beta$ factorisation of $w$.}
\end{figure}

\begin{proposition}\label{notu}
Let $u\in\Sigma^k$. Then $u\not\in\ScatFact_k(w)$ iff 
 $u[1] \in 
\alph(\beta_1) \setminus \alph(\alpha_1)$, $u[i] \in \alph(\beta_i)$, $u[i] u[i+1] 
\not\in\ScatFact_2(\beta_i \alpha_{i+1})$ for all 
$i\in[k-1]\backslash\{1\}$, and $u[k] \not\in \alph(\r(w))$.
\end{proposition}

\ifpaper
\else
\begin{proof}

Let  $w \in \NUniv_{\Sigma, m, k}$ with $\iota(w)=k-1$. 
%

Next, we want to characterise each length $k$ absent scattered factor of $w$. Consider 
$u\in\Sigma^k$ with $u \not\in \ScatFact_k(w)$. 
Suppose $u[1]\in\alph(\r(w^R))$.
By $\iota(w^R)=k-1$, we may choose $u[i]$ from $\ar_{k-i+1}(w^R)$ and get $u\in\ScatFact_k(w)$. 
Thus, $u[1] \not\in \alph(\r(w^R))$ but since $u[1]\in\alph(\ar_1(w))$ we have $u[1]\in 
\alph(\beta_1) \setminus \alph(\alpha_1)$.

If $u[1]u[2]\in \ScatFact_2(\ar_1(w))$, we could choose $u[i]\in\ar_{i-1}(w)$, for all 
$i\in[k-1]\backslash\{1,2\}$, and would get $u\in\ScatFact_k(w)$. Analogously to the argumentation 
for $u[1]\not\in\alph(\alpha_1)$, we get $u[1]u[2] \not\in\ScatFact_2(\alpha_1\beta_2)$. This concludes the 
induction basis.

As induction hypothesis, assume for one fixed $i\in[k-1]$ and all $\ell\in[i-1]$
\begin{enumerate}
\item $u[\ell] \in \alph(\beta_\ell)$,
\item $u[\ell]u[\ell+1] \not\in\ScatFact_2(s_{\ell} \alpha_{\ell+1})$ 
\end{enumerate}
Consider $u[i]$. Suppose that $u[i] \in\alph(p_{i})$.  Then, we have $u[i-1]u[i]$ is a scattered 
factor of $\ar_{k-i-1}(w^R)^R$. As there are exactly $k-i$ arches preceding 
$\alpha_{i}$ in $w^R$, we have $u\in\ScatFact_k(w)$, a contradiction. Thus, $u[i] \in\alph(\beta_{i})$.
With $\ar_{k-i}(w^R) = (\beta_i \alpha_{i+1})^R$, we have $u[i+1]\in\alph(\ar_i(w))$. Suppose $u[i] u[i+1] 
\in\ScatFact_2(\beta_i \alpha_{i+1})$, thus $u[1..i+1] \in\ScatFact_{i+1}(\ar_1(w) \cdots \ar_i(w) 
\alpha_{i+1})$ and as there are exactly $k-i-1$ arches preceding in $w^R$ it follows that 
$u\in\ScatFact_k(w)$, a contradiction. Thus, $u[i] u[i+1] \not\in\ScatFact_2(\beta_i \alpha_{i+1})$. This 
proves 1. and 2. for $u[i]$ and concludes the induction.

Additionally, $u[k] \not\in\alph(\r(w))$ as this would contradict 
$u \not\in\ScatFact_k(w)$ and, analogously, $u[k-1]u[k] 
\not\in\ScatFact_2(\beta_{k-1}\r(w))$.

So far, we haven proven that $u\not\in\ScatFact_k(w)$ then 
 $u[i] \in\alph(\beta_i)$, $u[i] u[i+1] \not\in\ScatFact_2(\beta_i 
\alpha_{i+1})$ for all $i\in [k-1]$, and $u[k] \not\in\alph(\r(w))$.

\medskip

Now, we prove the other direction, i.e., if the conditions hold for some $u\in\Sigma^k$, we 
have $u\not\in\ScatFact_k(w)$. Thus, consider $u\in\Sigma^k$ such that $u[1] \in 
\alph(\beta_1) \setminus \alph(\alpha_1)$, $u[i] \in \alph(\beta_i)$, $u[i] u[i+1] \not\in\ScatFact_2(\beta_i 
\alpha_{i+1})$ for all $i\in[k-1]\backslash\{1\}$, and $u[k] \not\in \alph(\r(w))$. Suppose that 
$u\in\ScatFact_k(w)$. Then, $u[1]$ occurs first in $\beta_1$, and the letters $u[i]$ can only be chosen 
first from $\ar_i(w)$ for every $i\in[k-1]$. But $u[k]\not\in\alph(\r(w))$ leads to a 
contradiction. 
Thus, $u\not\in\ScatFact_k(w)$.\qed
\end{proof}

\fi

Define $f_w:[|w|]\rightarrow[k]$ such that $f_w(i)=\ell$ iff $w$'s \nth{$i$} letter belongs to 
$\ar_{\ell}(w)$, for $i\in[|w|-|\r(w)|]$, and $f_w(i)=k$ otherwise. Moreover, define 
$g_{w,\ell}:\Sigma\rightarrow [|w|]$ by $g_{w,\ell}(\ta) =  \min\{ i\ |\ w[i] = \ta\ \land\ f_w(i) 
= \ell 
\}$ for all $\ell \in[k-1]$. Set $M_{w,1} = \alph(\beta_1) \setminus \alph(\alpha_1)$ and 
$M_{w,j} = (\alph(\beta_{i+1}) \setminus \alph(\beta_i[j'+1 .. |\beta_i|]\alpha_{i+1})) \cap 
\alph(\beta_i[1..j'])$ where $f_w(j) = i$, $j' = j - (\sum_{l=1}^{i} |\ar_l(w)| + |\alpha_i|)$ and $\alph(\beta_k) = \Sigma \setminus \alph(\r(w))$,
as well as $M_{w,1}^\prime = g_{w,1}(M_{w,1})$ and $M_{w,j}^\prime = g_{w,f(j)+1}(M_{w,j})$ for all $2 
\leq j < \max\{ m\ | f_w(m) < k \}$. Let $h_w(i) = \sum_{j \in M'_{w,i}} h_w(j)$ for all $i \in \{ 
\ell\mid f_w(\ell) < k-1 \}$ and $h_w(i) = |\Sigma \setminus \alph(\r(w))|$ otherwise.

\begin{remark}\label{modusinM}
Notice that by the definition of $\m(w)$ and Proposition~\ref{notu}, we have 
$\m(w)[i+1]\in M_{w,g_{w,f_w(i)}(\m(w)[i])}$ for all $i\in[k-2]$.
\end{remark}

\begin{proposition}\label{amount}
If $w\in\NUniv_{\Sigma,m,k}\cap\Univ_{\Sigma,k-1}$ then $m=h_w(1)$.
\end{proposition}

\ifpaper
\else
\begin{proof}
Choose a sequence of numbers $\mathcal{I}_u \in\N^{k-1}$ such that $\mathcal{I}_u[1] \in M_{w,1}'$ 
and $\mathcal{I}_u[i+1] \in M_{w,\mathcal{I}_u[i]}'$ for all $i\in[k-1]$.
Then for the word $u\in\Sigma^k$ such that $u[i] = w[\mathcal{I}_u[i]]$ for $i\in[k-1]$ and 
$u[k] \in \Sigma\setminus\alph(\r(w))$ we have
\begin{itemize}
\item $u[1] \in \alph(\beta_1) \setminus \alph(\alpha_1)$, 
\item $u[i] \in \beta_i$,
\item $u[i] u[i+1] \not\in\ScatFact_2(\beta_i \alpha_{i+1})$ for all $i\in[k-2]\backslash\{1\}$, and
\item $u[k] \not\in \alph(\r(w))$.
\end{itemize}
Thus, by Proposition~\ref{notu}, we have $u\not\in\ScatFact_k(w)$.
Then, calculating $h_w(1)$ recursively equals the number of possibilities to choose such sequences 
$\mathcal{I}_u$ and extend them with any letter $\ta_w \not\in\r(w)$. Each such sequence is 
associated to a different absent scattered factor $u$, i.e., $h(1)$ equals exactly the number $m$ 
of length $k$ absent scattered factors in $w$. \qed
\end{proof}

\fi

The following lemma shows that $u\in\Sigma^k$ is absent in $w,w'$ iff the sets of possible
candidates for positions in $\beta_i$ coincide for $w$ and $w'$ resp.

\begin{lemma}\label{umissing}
Let $w,w'\in\Univ_{\Sigma,k-1}$ with $\alph(\r(w'))\subseteq\alph(\r(w))$ and $u\in\Sigma^k$ with 
$u\not\in\ScatFact_k(w)$.
Choose $\mathcal{I}[1]\in M_{w,1}'$ and $\mathcal{I}[i+1]\in M_{w,\mathcal{I}[i]}'$
such that $u[1..k-1]=w[\mathcal{I}[1]]\cdots w[\mathcal{I}[k-1]]$. Then $u\not\in\ScatFact_k(w')$
iff there exist $\mathcal{I}'[1]\in M_{w',1}'$ and $\mathcal{I}'[i+1]\in M_{w',\mathcal{I}'[i]}'$
with $u[1..k-1]=w'[\mathcal{I}'[1]]\cdots w[\mathcal{I}'[k-1]]$ and $u[i]\in 
M_{w,\mathcal{I}[i]}\cap
M_{w',\mathcal{I}'[i]}$ for all $i\in[k-1]$.
\end{lemma}

\ifpaper
\else
\begin{proof}
First, consider $u\not\in\ScatFact_k(w')$. By the definition of the sets $M_{w',\cdot}'$, 
$M_{w,\cdot}$, and $M_{w',\cdot}$ the claim follows immediately.

The second direction follows by $\alph(\r(w'))\subseteq\alph(\r(w))$.\qed
\end{proof}

\fi

For $w,w'\in\Univ_{\Sigma,k-1}$ and $u\in\Sigma^k$, let $C(u,w,w')$
be the predicate of the iff-conditions for  $u\not\in\ScatFact_k(w)$.

\begin{theorem}\label{mequiv}
For all $w,w'\in\NUniv_{\Sigma,m,k}\cap\Univ_{\Sigma,k-1}$, we have $w\sim_k w'$ iff $C(u,w,w')$ 
and $C(u,w',w)$ for all $u\in\Sigma^k$. \end{theorem}

\ifpaper
\else
\begin{proof}
Follows directly by Lemma~\ref{umissing}.\qed
\end{proof}

\fi

Notice that $w\sim_k w'$ is equivalent to $M_{w,j}=M_{w',j'}$ for all $j,j'$ according to 
appropriate
sequences $\mathcal{I}$ and $\mathcal{I}'$ - illustrated in the following example.

\ifpaper
\else 
\begin{example}\label{example}
To give an example, consider the word $ w= \mathtt{(aabc) \cdot (bcca) \cdot b} \in \NUniv_{\Sigma, 4, 3}$. 
Applying Proposition~\ref{notu} results in the absent scattered factors $\mathtt{baa}, \mathtt{bac}, \mathtt{caa}, \mathtt{cac}$. Considering the appropriate factorisation 
in  $\alpha_{i}, \beta_{i}$ for $i \in [k-1]$, we get $\alpha_1 = \ta$, $\beta_1 = \mathtt{abc}$, $\alpha_2 = \mathtt{bc}$, $\beta_2 = \mathtt{ca}$  and $\alpha_3 = \mathtt{b}$.

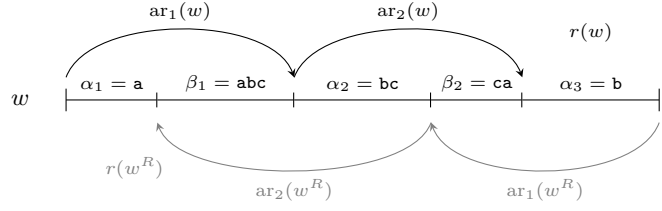
\begin{figure}
 \centering
\begin{tikzpicture}[scale=0.6,
blackbox/.style={rectangle, draw=black, thick, minimum height=0.6cm},
redbox/.style={rectangle, draw=red!75, 
			 thick, minimum height=0.5cm}
]

\coordinate (start) at (0,1);
\coordinate (end) at (13,1);

\coordinate (ar1) at (5,1);
\coordinate (ar2) at (10,1);
\coordinate (r) at (15,1);

\coordinate (rR) at (0,1);
\coordinate (ar2R) at (2,1);
\coordinate (ar1R) at (8,1);

\node[black] (w) at (-1,1) {$w$};

\draw[-,black] ($ (start) - (0,0.25) $) -- ($ (start) + (0,0.25) $) ;
\draw[-,black] ($ (end) - (0,0.25) $) -- ($ (end) + (0,0.25) $) ;

\draw[black, -stealth] ($ (start)+(0,0.5) $) to[out=70, in=110, distance=1.5cm] node[above, black, font=\scriptsize] {$\ar_1(w)$} ($ (ar1)+(0,0.5) $);
\draw[black, -stealth] ($ (ar1)+(0,0.5) $) to[out=70, in=110, distance=1.5cm] node[above, black, font=\scriptsize] {$\ar_2(w)$} ($ (ar2)+(0,0.5) $);
\node[black, font=\scriptsize] (r) at ($ (end) + (-1.5,1.5) $) {$r(w)$};

\draw[gray, -stealth] ($ (end)+(0,-0.5) $) to[out=-110,in=-70, distance=1.5cm] node[below, gray, font=\scriptsize] {$\ar_1(w^R)$}  ($ (ar1R)+(0,-0.5) $);
\draw[gray, -stealth] ($ (ar1R)+(0,-0.5) $) to[out=-110,in=-70, distance=1.5cm] node[below, gray, font=\scriptsize] {$\ar_2(w^R)$} ($ (ar2R)+(0,-0.5) $);
\node[gray, ->, font=\scriptsize] (rR) at ($ (start) + (1.5,-1.5) $) {$r(w^R)$};

\draw[-|,black] (start) -- node[above,font=\scriptsize] {$\alpha_1 = \ta$} (ar2R);
\draw[-|,black] (ar2R) -- node[above,font=\scriptsize] {$\beta_1=\ta\tb\tc$} (ar1);
\draw[-|,black] (ar1) -- node[above,font=\scriptsize] {$\alpha_2=\tb\tc$} (ar1R);
\draw[-|,black] (ar1R) -- node[above,font=\scriptsize] {$\beta_2=\tc\ta$} (ar2);
\draw[-,black] (ar2) -- node[above,font=\scriptsize] {$\alpha_3=\tb$} (end);



\end{tikzpicture}
\caption{Factorisation of $w = \mathtt{aabcbccab}$}
\end{figure}

Now, we want to calculate $h_w(1)$ as in Proposition~\ref{amount}. Thus, we need to consider $f_w$ 
first and have	$f_w(i) = 1 \text{ for } i \in [4]$, $f_w(i) = 2 \text{ for } i \in 
[8]\backslash[4]$ and $f_w(9) = 3$.
Now, $g_{w,l}(\ta)$ defines the index of the leftmost occurrence of $\ta$ in the \nth{$\ell$} arch. Here we give an example for the leftmost occurrence of $\ta$ in the first arch, 
described by $g_{w,1}(\ta) = \min\{i \mid w[i] = \ta \land f_w(i) = 1\} = \min\{1,2\} = 1$, and in the second arch respectively, i.e., $g_{w,2}(\ta) = \min\{i \mid w[i] = \ta \land f_w(i) = 2\} = \min\{8\} = 8$. 
By definition we have
\begin{align*}
	M_{w,1} &= \alph(\beta_1) \setminus \alph(\alpha_1) = \alph(\mathtt{abc}) \setminus \alph(\mathtt{a}) = \{\tb,\tc\},\\
	M_{w,3} &= (\alph(\beta_{1+1}) \setminus \alph(\beta_1[2+1 .. |\beta_1|]\alpha_{1+1})) \cap \alph(\beta_1[1..2])\\
			&= (\{\tc, \ta\} \setminus \{\tc, \tb\}) \cap \{\ta, \tb\} = \{\ta\},\\
	M_{w,4} &= (\alph(\beta_{1+1}) \setminus \alph(\beta_1[3+1 .. |\beta_1|]\alpha_{1+1})) \cap \alph(\beta_1[1..3])\\
			&= (\{\tc, \ta\} \setminus \{\tc, \tb\}) \cap \{\ta, \tb\, \tc \} = \{\ta\},\\
	M_{w,8} &= (\alph(\beta_{2+1}) \setminus \alph(\beta_2[2+1 .. |\beta_1|]\alpha_{2+1})) \cap \alph(\beta_2[1..2])\\
			&= (\{\ta, \tc\} \setminus \{\tb\}) \cap \{\tc, \ta\} = \{\ta, \tc\}.		
\end{align*}

Further, we get
\begin{align*}
	M'_{w,1} &= g_{w,1}(M_{w,1}) = g_{w,1}(\{\tb,\tc\}) = \{3,4\},\\
	M'_{w,3} &= g_{w,f(3) +1}(M_{w,3}) = g_{w,2}(\{\ta\}) = \{8\},\\
	M'_{w,4} &= g_{w,f(4) +1}(M_{w,4}) = g_{w,2}(\{\ta\}) = \{8\}.
\end{align*}

Notice that $M'_{w,j}$ for all $j < 4$ is not defined since $f_w(j) > 1$, thus $g_{w,f(j)+1}$ is not defined.
Now, it is easy to see how the sequences $\mathcal{I}_u$ belong to the absent scattered factors of $w$. With $\mathcal{I}_u[1] \in M'_{w,1}$ and $\mathcal{I}_u[2] \in M'_{w,\mathcal{I}_u[1]}$, the possible sequences are $(3,8)$ and $(4,8)$. Since $w[3] = \tb$, $w[4] = \tc$ and $w[8] = \ta$ all absent scattered factors of $w$ have either one of them as prefix and end in one of the letters $\ta$ or $\tc$ (missing in $\r(w)$).

To determine $h_w(1)$, we have with $h_w(8) = |\Sigma| - |\alph(\r(w))| = 2$
\begin{align*}
h_w(1) & = \sum_{j \in M'_{w,1}}  h_w(j)\\ 
&= h_w(3) + h_w(4)\\
&= \sum_{j \in M'_{w,3}}
h_w(j) + \sum_{j \in M'_{w,4}} h_w(j)  \\
		&= h_w(8) + h_w(8)\\
		&= 4.
\end{align*}

Moreover, we have $\mathtt{(aabc) \cdot (bcca) \cdot b} \sim_k \mathtt{(aabbc) \cdot (bccca) \cdot b}$ since the 
letters occurring in $\alpha'_i$, $\beta'_i$ of the factorisation of $\mathtt{abbcbcccab}$ for $i 
\in [k-1]$ are pairwise equal to those in $w$.
\begin{figure}
	\centering
\begin{tikzpicture}[scale=0.6,
blackbox/.style={rectangle, draw=black, thick, minimum height=0.6cm},
redbox/.style={rectangle, draw=red!75, 
			 thick, minimum height=0.5cm}
]

\coordinate (start) at (0,1);
\coordinate (end) at (13,1);

\coordinate (ar1) at (5,1);
\coordinate (ar2) at (10,1);
\coordinate (r) at (15,1);

\coordinate (rR) at (0,1);
\coordinate (ar2R) at (2,1);
\coordinate (ar1R) at (8,1);

\node[black] (w) at (-1,1) {$w'$};

\draw[-,black] ($ (start) - (0,0.25) $) -- ($ (start) + (0,0.25) $) ;
\draw[-,black] ($ (end) - (0,0.25) $) -- ($ (end) + (0,0.25) $) ;

\draw[black, -stealth] ($ (start)+(0,0.5) $) to[out=70, in=110, distance=1.5cm] node[above, black, font=\scriptsize] {$\ar_1(w')$} ($ (ar1)+(0,0.5) $);
\draw[black, -stealth] ($ (ar1)+(0,0.5) $) to[out=70, in=110, distance=1.5cm] node[above, black, font=\scriptsize] {$\ar_2(w')$} ($ (ar2)+(0,0.5) $);
\node[black, font=\scriptsize] (r) at ($ (end) + (-1.5,1.5) $) {$r(w')$};

\draw[gray, -stealth] ($ (end)+(0,-0.5) $) to[out=-110,in=-70, distance=1.5cm] node[below, gray, font=\scriptsize] {$\ar_1(w'^R)$}  ($ (ar1R)+(0,-0.5) $);
\draw[gray, -stealth] ($ (ar1R)+(0,-0.5) $) to[out=-110,in=-70, distance=1.5cm] node[below, gray, font=\scriptsize] {$\ar_2(w'^R)$} ($ (ar2R)+(0,-0.5) $);
\node[gray, ->, font=\scriptsize] (rR) at ($ (start) + (1.5,-1.5) $) {$r(w'^R)$};

\draw[-|,black] (start) -- node[above,font=\scriptsize] {$\alpha'_1 = \ta$} (ar2R);
\draw[-|,black] (ar2R) -- node[above,font=\scriptsize] {$\beta'_1=\ta\tb\tb\tc$} (ar1);
\draw[-|,black] (ar1) -- node[above,font=\scriptsize] {$\alpha'_2=\tb\tc$} (ar1R);
\draw[-|,black] (ar1R) -- node[above,font=\scriptsize] {$\beta'_2=\tc\ta$} (ar2);
\draw[-,black] (ar2) -- node[above,font=\scriptsize] {$\alpha'_3=\tb$} (end);



\end{tikzpicture}
\caption{Factorisation of $w' = \mathtt{aabbcbccab}$}
\end{figure}

Similarly, we have $\mathtt{aabcbccab} \not\sim_k  \mathtt{aabcbcab} \in \NUniv_{\Sigma, 5, 3}$, 
assuming a given factorisation into $\alpha''_i$ and $\beta''_i$ for $i \in [k-1]$, since $\alpha_2 
\neq \alpha''_2$ as illustrated in Figure~\ref{word3}. Thus, $\mathtt{ccc}$ is absent as well.

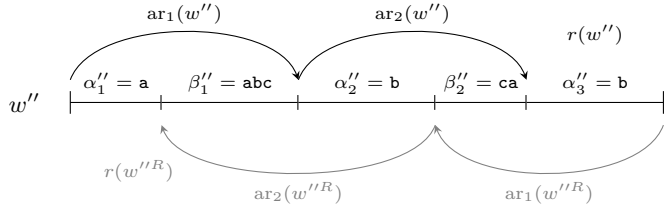
\begin{figure}
	\centering
\begin{tikzpicture}[scale=0.6,
blackbox/.style={rectangle, draw=black, thick, minimum height=0.6cm},
redbox/.style={rectangle, draw=red!75, 
			 thick, minimum height=0.5cm}
]

\coordinate (start) at (0,1);
\coordinate (end) at (13,1);

\coordinate (ar1) at (5,1);
\coordinate (ar2) at (10,1);
\coordinate (r) at (15,1);

\coordinate (rR) at (0,1);
\coordinate (ar2R) at (2,1);
\coordinate (ar1R) at (8,1);

\node[black] (w) at (-1,1) {$w''$};

\draw[-,black] ($ (start) - (0,0.25) $) -- ($ (start) + (0,0.25) $) ;
\draw[-,black] ($ (end) - (0,0.25) $) -- ($ (end) + (0,0.25) $) ;

\draw[black, -stealth] ($ (start)+(0,0.5) $) to[out=70, in=110, distance=1.5cm] node[above, black, font=\scriptsize] {$\ar_1(w'')$} ($ (ar1)+(0,0.5) $);
\draw[black, -stealth] ($ (ar1)+(0,0.5) $) to[out=70, in=110, distance=1.5cm] node[above, black, font=\scriptsize] {$\ar_2(w'')$} ($ (ar2)+(0,0.5) $);
\node[black, font=\scriptsize] (r) at ($ (end) + (-1.5,1.5) $) {$r(w'')$};

\draw[gray, -stealth] ($ (end)+(0,-0.5) $) to[out=-110,in=-70, distance=1.5cm] node[below, gray, font=\scriptsize] {$\ar_1(w''^R)$}  ($ (ar1R)+(0,-0.5) $);
\draw[gray, -stealth] ($ (ar1R)+(0,-0.5) $) to[out=-110,in=-70, distance=1.5cm] node[below, gray, font=\scriptsize] {$\ar_2(w''^R)$} ($ (ar2R)+(0,-0.5) $);
\node[gray, ->, font=\scriptsize] (rR) at ($ (start) + (1.5,-1.5) $) {$r(w''^R)$};

\draw[-|,black] (start) -- node[above,font=\scriptsize] {$\alpha''_1 = \ta$} (ar2R);
\draw[-|,black] (ar2R) -- node[above,font=\scriptsize] {$\beta''_1=\ta\tb\tc$} (ar1);
\draw[-|,black] (ar1) -- node[above,font=\scriptsize] {$\alpha''_2=\tb$} (ar1R);
\draw[-|,black] (ar1R) -- node[above,font=\scriptsize] {$\beta''_2=\tc\ta$} (ar2);
\draw[-,black] (ar2) -- node[above,font=\scriptsize] {$\alpha''_3=\tb$} (end);



\end{tikzpicture}
\caption{Factorisation of $w'' = \mathtt{aabcbcab}$}
	\label{word3}
\end{figure}
\end{example}

\fi

In this section we showed for some $m$ how $\NUniv_{\Sigma,m,k}$ looks like and determined
$m$ for $w\in\NUniv_{\Sigma,m,k}\cap\Univ_{\Sigma,k-1}$ as well as $[w]_{\sim_{k}}$.

\section{Conclusion}\label{conc}
In this work, we pursued the approach to partition $\Sigma^{\ast}$ w.r.t. the number of absent 
scattered factors of a given length $k$. This lead to the notion of $m$-nearly $k$-universal words, 
which are words where exactly $m$ scattered factors of length $k$ are absent. We haven chosen this
perspective to investigate the index of the Simon congruence $\sim_k$ and indeed we were able to
fully characterise $1$-nearly $k$-universal words and give the index as well as a characterisation 
of  $\sim_k$ restricted to this subclass. Moreover, we gave some insights for $m>1$, especially for 
$m\in\{2,|\Sigma|^k-1,|\Sigma|^k-2\}$ (notice that $m=0$ is fully investigated in 
\cite{barker2020scattered}). 
Additionally in Section~\ref{mnearly}, we followed the idea from \cite{kosche2021absent} from a 
combinatorial point of view, showing for instance that letters have the same dist-value in 
\cite{kosche2021absent} iff they are in the same arch of $w^R$. By this approach we showed that
$m$ can be determined recursively for $w$ with $\iota(w)=k-1$ by investigating the overlaps of the 
arches from $w$ and $w^R$. Moreover, we proved when to words $w_1,w_2$ with 
$\iota(w_1)=\iota(w_2)=k-1$ fulfil $w_1\sim_k w_2$.

Unfortunately, we were not able to give a full characterisation of $\NUniv_{\Sigma,m,k}$ for 
arbitrary $m$. A first step could be to determine $\iota(w)$ for $w\in\NUniv_{\Sigma,m,k}$.
We conjecture that chosing $i\in[\sigma^k]$ such that $\sigma^i\leq m\leq\sigma^{i+1}-1$ leads to  
$k - \lfloor \frac{m}{\sigma^i + 1} \rfloor - 1 \leq \iota(w) \leq k - \lfloor 
\frac{m}{\sigma^i + 1} \rfloor$. A subpartition of $\NUniv_{\Sigma,m,k}$ depending on $\iota(w)$
(as introduced in \cite{kosche2021absent} and used in Section~\ref{mnearly}) could prove useful.

\bibliography{refs}


\end{document}